\documentclass{amsart}

\usepackage{graphicx}
\usepackage{amscd}
\usepackage{amsmath}
\usepackage{amsfonts}
\usepackage{amssymb}
\usepackage{amsxtra}
\usepackage{xspace}
\usepackage{array}

\theoremstyle{plain}
\newtheorem{theorem}{Theorem}
\newtheorem{lemma}{Lemma}
\newtheorem{corollary}{Corollary}
\newtheorem{proposition}{Proposition}

\newtheorem{claim}{Claim}
\newtheorem{example}{Example}
\newtheorem{remark}{Remark}

\theoremstyle{definition}
\newtheorem{definition}{Definition}

\numberwithin{equation}{section}

\newcommand{\aip}[2][(p)]{a_{i{#1}}(X_{#2})} 
\newcommand{\aipb}[2][(p)]{a_{\bi{#1}}(X_{#2})} 
\newcommand{\nob}[1][\gamma]{|#1|_{\xp}}
\newcommand{\noba}[1][\gamma]{|#1|}
\newcommand{\tx}[1][]{\widetilde{X}_{#1}}

\newcommand{\Nx}[1][]{N(#1)} 

\newcommand{\kv}{\cite{KV}\xspace}
\newcommand{\yk}{Y_k}

\newcommand{\om}{\omega}

\newcommand{\omi}[1][i]{\omega_{#1}^{(n)}}
\newcommand{\ombi}[1][i]{\bar{\omega}_{#1}^{(n)}}
\newcommand{\tom}{\tilde{\omega}}
\newcommand{\ri}{\ensuremath{\; i=1,2}}
\newcommand{\bri}[1][k]{\widehat{Y}_{#1}}
\newcommand{\bi}{\,\bar{\imath}}
\newcommand{\N}{\mathbb{Z}^{\geq 0}}
\newcommand{\nk}[1][k]{N(Z_{#1})}
\newcommand{\D}[1][i]{\Delta_{#1}\,}
\newcommand{\dyn}{Dynkin diagram\xspace}
\newcommand{\dyns}{Dynkin diagrams\xspace}
\newcommand{\ep}{ extensible pair\xspace}
\newcommand{\al}{\alpha}
\newcommand{\ali}[1][i]{\alpha_{#1}^{(n)}}
\newcommand{\albi}[1][i]{\bar{\alpha}_{#1}^{(n)}}
\newcommand{\tal}{\tilde{\alpha}}
\newcommand{\alch}{\Check{\alpha}}
\newcommand{\integers}{\mathbb{Z}}
\newcommand{\rationals}{\mathbb{Q}}

\newcommand{\be}{\begin{enumerate}}
\newcommand{\ee}{\end{enumerate}}
\newcommand{\beq}{\begin{equation}}
\newcommand{\eeq}{\end{equation}}
\newcommand{\bprop}{\begin{proposition}}
\newcommand{\eprop}{\end{proposition}}

\newcommand{\omuv}{\om_u - \om_v}
\newcommand{\omvu}{\om_v - \om_u}

\newcommand{\kma}{\mathfrak{g}}
\newcommand{\csa}{\mathfrak{h}}
\newcommand{\xp}{\ensuremath{(X_1, X_2)\;}}
\newcommand{\hht}[1][]{{\mathop{\mathrm{ht}}\nolimits}_{#1}}

\newcommand{\clmn}[1][k]{c_{\lambda\mu}^{\, \nu}(#1)}
\newcommand{\clmng}[1][k]{c_{\lambda\mu\nu}^{\; \gamma}(#1)}
\renewcommand{\L}{\lambda}    
\newcommand{\Lk}{\lambda\rk}    
\newcommand{\rk}[1][k]{^{(#1)}}
\newcommand{\mk}{\mu\rk}    

\newcommand{\W}{W}    
\newcommand{\R}{R}

\newcommand{\Pk}{\Pi\rk}
\newcommand{\Sk}[1][k]{\Sigma_{#1}}

\newcommand{\plk}{\pi_{\lambda\rk}}

\newcommand{\pk}[1][\beta]{\mathcal{P}_k(\lambda, #1)}
\newcommand{\pkp}[1][\beta]{\mathcal{P}^+_k(\lambda, #1)}
\newcommand{\pkm}[1][\beta]{\mathcal{P}_k^{(\mu)}(\lambda, #1)}
\newcommand{\pko}[1][\beta]{\mathcal{P}^0_k(\lambda, #1)}

\newcommand{\pkppr}[1][\beta]{\mathcal{P}^+_{k\pr}(\lambda, #1)}
\newcommand{\pkmpr}[1][\beta]{\mathcal{P}_{k\pr}^{(\mu)}(\lambda, #1)}

\newcommand{\ho}{\mathcal{H}_1}
\newcommand{\hop}{\mathcal{H}^+_1}
\newcommand{\Ht}{\mathcal{H}_2}
\newcommand{\Htp}{\mathcal{H}^+_2}

\newcommand{\ls}{{\mathop{\mathrm{ls}}\nolimits}} 
\newcommand{\rs}{{\mathop{\mathrm{rs}}\nolimits}} 

\newcommand{\xip}{x_{i(p)}}
\newcommand{\yipb}{y_{\bi(p)}}

\newcommand{\pzk}[1][]{P^{#1}(Z_k)}
\newcommand{\qzk}[1][]{Q^{#1}(Z_k)}

\newcommand{\g}{\gamma}
\newcommand{\gk}[1][k]{\gamma^{(#1)}}

\newcommand{\blb}[1][\beta]{b_{\lambda #1}}

\newcommand{\T}{{\mathop{\mathrm{ht}}\nolimits}}

\newcommand{\po}{\succcurlyeq}
\newcommand{\mcl}{\ensuremath{\mathcal{R}\;}}
\newcommand{\mclh}{\ensuremath{\widehat{\mathcal{R}}\;}}
\newcommand{\rxx}{\mathcal{R}(X_1|X_2)}
\newcommand{\rxxa}{\mathcal{R}(X|A_1)}

\newcommand{\bxo}[1][\lambda]{(#1)_{X_1}}
\newcommand{\bxt}[1][\lambda]{(#1)_{X_2}}
\newcommand{\pr}{^\prime}
\newcommand{\ppr}{^{\prime\prime}}

\newcommand{\complex}{\mathbb{C}}
\newcommand{\supp}{\ell}
\newcommand{\dep}{{\mathop{\mathrm{dep}}\nolimits}}

\begin{document}

\title[Dynkin diagram sequences]{Dynkin diagram sequences and stabilization phenomena}
\author{Sankaran Viswanath}
\address{Department of Mathematics\\
University of California\\
Davis, CA 95616, USA}
\email{svis@math.ucdavis.edu}
\subjclass{17B67}
\keywords{Tensor product multiplicity, branching multiplicity,
  Littelmann path model, partial order on dominant weights}

\begin{abstract}
We continue the study of stabilization phenomena for Dynkin diagram
sequences initiated in the earlier work of Kleber and the present
author. We consider a more general class of sequences than that of
this earlier work, and isolate a condition on the weights that gives
stabilization of tensor product and branching multiplicities. We show
that all the results of the previous article can be naturally
generalized to this setting.
We also prove some properties of the partially ordered set of dominant
weights of indefinite Kac-Moody algebras, and use this to give a more
concrete definition of a stable representation ring. Finally, we
consider the classical sequences $B_n, C_n, D_n$ that fall outside the
purview of the earlier work, and work out some easy-to-describe
 conditions on the weights which imply stabilization.
\end{abstract}
\maketitle


\section{Introduction}
In this article, we will consider
sequences of Dynkin diagrams and study the behavior of representations
of the associated Kac-Moody algebras. The sequences of Dynkin diagrams
considered are of the form 
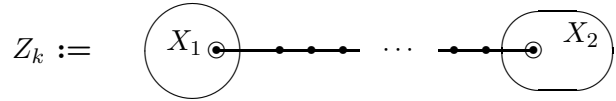
\setlength{\unitlength}{12pt}
$$\begin{picture}(5,3)(2,-1)
\put(-5,0){\makebox(0,0){\bf \large $Z_k$ := }}
\put(-1,0.25){\makebox(0,0){\bf \large $X_1$}}
\put(-.75,0){\circle{3}}
\put(0,0){\circle*{.25}}
\put(0,0){\circle{.5}}
\put(2,0){\circle*{.25}}
\put(3,0){\circle*{.25}}
\put(4,0){\circle*{.25}}
\put(7.5,0){\circle*{.25}}
\put(8.5,0){\circle*{.25}}
\put(10,0){\circle*{.25}}
\put(10,0){\circle{.5}}
\put(10.75,0){\oval(3.5,2.5)}
\put(11.5,0.5){\makebox(0,0){\bf \large $X_2$}}
\put(0,0){\line(1,0){4.5}}
\put(10,0){\line(-1,0){3}}
\put(5.3,0){\makebox(1,0){$\cdots$}}

\end{picture}$$
where $X_1$
and $X_2$ are fixed  Dynkin diagrams and the string of
intermediate nodes has length $k$. The article \kv considered the
Dynkin diagram sequences $Z_k$ arising in the special case 
when $X_2=A_1$ (the diagram with just a single node).

For most choices of $X_1, X_2$, the associated Kac-Moody algebra
$\kma(Z_k)$ is infinite dimensional, non-affine (i.e  of {\em
  indefinite} type) and very little is known about such Lie algebras
in general. So, rather than study representations of 
the individual $\kma(Z_k)$, we study them in the limit as $k \to \infty$. 

We consider dominant integral weights $\lambda, \mu, \nu$ which are
{\em supported} on the two ends of the Dynkin diagram of $Z_k$. In
\kv, the primary object of interest was the multiplicity $\clmn$ of 
the irreducible highest weight representation $L(\nu\rk)$ of $\kma(Z_k)$ 
in the tensor product $L(\lambda\rk) \otimes L(\mu\rk)$. Specifically
under some additional 
conditions, it was shown that these multiplicities become constant  ({\em
  stabilize}) for large $k$; this generalizes the classical
stabilization results for the $A_n$ diagram which more or less follow
from the Littlewood-Richardson rule.

In this article, we consider  the tensor product multiplicity $\clmn$ as well as the
branching multiplicities $\blb(k)$. The question we ask is this :

\noindent
Under what conditions on the weights $\L, \mu, \nu, \beta$ 
do $\clmn$ and $\blb(k)$ stabilize ?  

In section \ref{thisone} we provide a sufficient condition that
ensures stabilization.
Having formulated this criterion for stabilization, we follow two distinct
threads:

First, we recall that the goal of the previous article \kv was slightly
different; it sought to find conditions on the {\em diagram} $X_1$
($X_2= A_1$ there) which would ensure stabilization of $\clmn$ for
{\em all} weights $\L, \mu, \nu$. The condition on $X_1$ which made
this work  was called {\em extensibility}. In sections 
\ref{extpair} and \ref{h2plus}, we generalize this approach
and define a notion of extensibility for pairs $\xp$. This new notion
is strong enough to ensure stabilization for all $\L, \mu, \nu, \beta$
while still admitting enough interesting examples of diagrams. In
particular, this notion subsumes the earlier notion of \kv.

It was also shown in \kv that one could use the stable values of the
$\clmn$ to define an operation $*$, which mimics the limit as $k \to
\infty$ of the tensor product. A very surprising fact discovered there
was the associativity of $*$. A notion of a stable representation ring
was formulated as a consequence. We derive all these results for
extensible pairs \xp and define a more concrete, modified version of
the stable representation ring in this case. This appears in section
\ref{stabrepring}. 

As our second thread, we turn to the classical sequences of \dyns $B_n, C_n, D_n$. These were notable
   exceptions to the extensibility condition of \kv. So, while 
   nice stabilization results hold for the $A_n$, nothing much could be said about
   these other classical types. We remedy this in section \ref{BCD}. For
   these types, our methods do not imply stabilization for all choices of $\L,
   \mu, \nu, \beta$, but we work out some easy to describe conditions
   on the weights under which they do.

While this article was in preparation, Webster \cite{ben} has shown
that a more general version of our stabilization result (Theorem \ref{mainthm1}) can
be proved using quiver varieties and their connections with
representations of Kac-Moody algebras. The results and formulation in \cite{ben}
and in section \ref{thisone} of the present article overlap substantially.
Webster's approach also proves a `polynomiality of weight multiplicities' result
for these kinds of \dyn sequences.

\noindent
{\bf Acknowledgements}: I'd like to thank Michael Kleber for his
constant encouragement and many valuable suggestions on an earlier
draft of this article. I'd also like to thank Ben Webster for 
explaining many aspects of his approach to me. 

\section{A criterion for stabilization}\label{thisone}

\subsection{Notations}\label{nota}
We begin with some notations concerning Kac-Moody algebras and Dynkin
diagrams. Let $X$ be a \dyn (in the sense of 
\cite[Chapter 4]{Kac}) associated to a symmetrizable generalized Cartan matrix $C(X)$.
Using the data of $X$, one
constructs $\kma(X)$, the Kac-Moody algebra 
associated to $X$.  Let $\csa(X)$ denote the Cartan subalgebra of
 $\kma(X)$, with dual $\csa^*(X)$. We let $\Nx[X]$ denote  the set of
nodes of $X$ and $\det X:=\det C(X)$. By abuse of notation we will usually
write $p \in X$ to mean that $p \in \Nx[X]$.
For each $p \in X$, let
 $\alch_p, \al_p, \om_p$ respectively denote the
simple coroot, the simple root and the fundamental weight
 corresponding to $p$. 
So for example $\om_p(\alch_q) = \delta_{pq}$ for $p,q \in X$.

We let $Q(X), P(X)$ be the root and weight lattices of $\kma(X)$. The set 
\{$\al_p$, $p \in X\}$ forms a $\integers$ basis of $Q(X)$, 
and when $\det X \neq 0$ the set  $\{\om_p$, $p\in X\}$ forms a $\integers$ basis of $P(X)$.
Let $Q^+(X)$ be the set of $\integers^{\geq 0}$ linear combinations of
the $\al_p$ and $P^+(X)$ denote the set of  $\integers^{\geq 0}$
linear combinations of the $\om_p$.

When $\det X=0$, the fundamental weights are not uniquely defined and
we pick one of the possible $\om_p$'s for each node $p$. We will
usually not run into diagrams with $\det X=0$ (see section
\ref{finalrem} for an exception).

\subsection{Dynkin diagram sequences}
In this article, we will be interested in {\em sequences}
 of Dynkin diagrams; for example, we have the classical Dynkin diagram 
sequence $A_k$, $k \geq 1$ 
$$
\setlength{\unitlength}{15pt}
\begin{picture}(7,1)(1,-.5)
\put(0,0){\circle*{.25}}
\put(1,0){\circle*{.25}}
\put(2,0){\circle*{.25}}
\put(5,0){\circle*{.25}}
\put(0,0){\line(1,0){3}}
\put(5,0){\line(-1,0){1}}
\put(3.05,0){\makebox(1,0){...}}
\put(0,-.5){\makebox(0,0){\tiny 1}}
\put(1,-.5){\makebox(0,0){\tiny 2}}
\put(2,-.5){\makebox(0,0){\tiny 3}}
\put(5,-.5){\makebox(0,0){\tiny $k$}}
\end{picture}$$
We will let $A_0$ denote the empty diagram.

By a {\em marked Dynkin diagram}, we will mean the data $(X, \xi)$
where $X$ is a \dyn and $\xi \in X$ is a distinguished node.

Given a marked \dyn $(X,\xi)$ and an integer $m \geq 0$, one can
construct the \dyn $X(m)$ obtained by ``attaching'' 
the diagram $A_m \,(m \geq 0)$ to the node $\xi$ as follows:
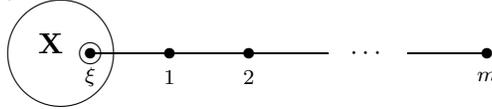
\setlength{\unitlength}{15pt}
$$\begin{picture}(4,2)(2,-1)
\put(-1,0.25){\makebox(0,0){\bf \large X}}
\put(-.75,0){\circle{3}}
\put(0,0){\circle*{.25}}
\put(0,0){\circle{.5}}
\put(2,0){\circle*{.25}}
\put(4,0){\circle*{.25}}
\put(10,0){\circle*{.25}}
\put(0,0){\line(1,0){6}}
\put(10,0){\line(-1,0){2}}
\put(6.5,0){\makebox(1,0){$\cdots$}}
\put(0,-.6){\makebox(0,0){\footnotesize $\xi$}}
\put(2,-.6){\makebox(0,0){\footnotesize $1$}}
\put(4,-.6){\makebox(0,0){\footnotesize  $2$}}
\put(10,-.6){\makebox(0,0){\footnotesize  $m$}}
\end{picture}$$

In the notation of \kv, this is $X_{d+m}$ where $d$ is the number of
nodes in $X$. We let $\xi(m)$ denote the end node (labeled $m$ in the
figure) and consider $(X(m), \xi(m))$ as a marked \dyn.

One can consider the sequence of symmetrizable Kac-Moody algebras
$\kma(X(m))$ associated with the $X(m)$. For most choices of $X$ and
for most values of $m$,
these turn out to be  infinite dimensional non affine Kac-Moody
algebras (i.e of  {\em indefinite type}), but one can still study
their integrable  highest weight representations. The objective of \kv
was to study how multiplicities in tensor product decompositions of
such representations of $\kma(X(m))$ change with $m$.

\subsection{Pairs}
We now consider a broader class of sequences of \dyns; these will be
obtained starting with {\em pairs} \xp of marked \dyns, rather than
with a single diagram $X$.

For $i=1,2$ let $(X_i, \xi_i)$  be given marked Dynkin diagrams 
such that $C(X_i)$ are symmetrizable. 
For each $k \geq 1$ we form a Dynkin
diagram $Z_k = Z_k ( X_1, X_2)$ by taking the diagram $A_k$ and
attaching its two ends to $\xi_1$ and $\xi_2$ as shown in figure :
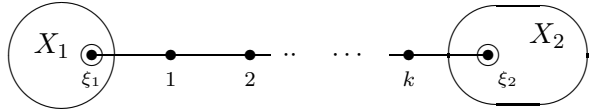
\setlength{\unitlength}{15pt}
$$\begin{picture}(5,3)(2,-1.5)\label{figu-orig}
\put(-1,0.25){\makebox(0,0){\bf \large $X_1$}}
\put(-.75,0){\circle{3}}
\put(0,0){\circle*{.25}}
\put(0,0){\circle{.5}}
\put(2,0){\circle*{.25}}
\put(4,0){\circle*{.25}}
\put(8,0){\circle*{.25}}
\put(10,0){\circle*{.25}}
\put(10,0){\circle{.5}}
\put(10.75,0){\oval(3.5,2.5)}
\put(11.5,0.5){\makebox(0,0){\bf \large $X_2$}}
\put(0,0){\line(1,0){4.5}}
\put(10,0){\line(-1,0){2.5}}
\put(6,0){\makebox(1,0){$\cdots$}}
\put(4.75,0){\makebox(0.5,0){$\cdot \cdot$}}
\put(0,-.6){\makebox(0,0){\scriptsize $\xi_1$}}
\put(2,-.6){\makebox(0,0){\footnotesize $1$}}
\put(4,-.6){\makebox(0,0){\footnotesize  $2$}}
\put(8,-.6){\makebox(0,0){\footnotesize  $k$}}
\put(10.5,-.6){\makebox(0,0){\scriptsize  $\xi_2$}}
\end{picture}$$
The figure doesn't show the rest of the nodes of $X_1$ and
$X_2$. The matrix $C(Z_k)$ is clearly symmetrizable. The associated
Kac-Moody algebra will be denoted $\kma(Z_k)$. 
If $X_2 = A_1$, the diagram with a single node, this construction
coincides with the earlier one; $Z_k = X_1(k+1)$.

In the rest of section \ref{thisone}, we will often need to refer to the following subdiagrams of the $Z_k$'s.
\be
\item 
$X_{12} := X_1 \cup X_2$. 
\item For $l \geq 0$, we  identify $X_1(l)$ with 
the subdiagram of $Z_k$ formed by $X_1$ and the  
intermediate nodes labeled $1,\cdots, l$ in the above figure.

\item Similarly given $r \geq 0$, identify $X_2(r)$ with the
subdiagram formed by $X_2$ and the intermediate nodes labeled 
$k, k-1,\cdots, k-r+1$. 

\vspace{.2cm}
Each of these is a subset of $Z_k$ for all but finitely many values of
$k$. When these subsets are encountered, the specific value of $k$
being used will be clear from context.    

\vspace{.1cm}
\item 
With $l, r$ as above, let $\yk(l,r) \subset Z_k$ denote
the subdiagram isomorphic to $A_{k-l-r}$ formed by the
nodes numbered $l+1, l+2, \cdots, k-r$.
So for instance, $\yk(0,0)$ is the diagram $A_k$ in the middle. 
\ee

\subsection{Representations of $\kma(Z_k)$}\label{twofour}
We now consider integrable highest weight representations of the
Kac-Moody algebra $\kma(Z_k)$. These are indexed by dominant integral
weights of $\kma(Z_k)$. Our immediate interest will be in dominant
integral weights which are ``supported'' on $X_1 \cup X_2$.
To make this more precise :

Let $\W\xp$ be the set of all functions $f: N(X_1) \cup N(X_2) \to
\integers$ and $\W^{+}\xp :=\{f \in \W\xp: image(f) \subset \N\}$. Similarly for each
$s \in \integers$, let $\R_s\xp$ be the set of all functions 
  $f: N(X_1) \cup N(X_2) \to\integers$ such that $f(\xi_1) = f(\xi_2)
= s$ and  $\R_s^+\xp:=\{f \in \R_s\xp: image(f) \subset \N\}$. We will let the
elements of $\W\xp$ and $\R_s\xp$ define elements of the weight and root
lattices of $\kma(Z_k)$ as follows: given $\lambda \in \W\xp$, define $\Lk
\in P(Z_k)$ by
$$ \Lk := \sum_{p \in X_{12} \subset Z_k} \L(p) \, \om_p$$
We define $P_0(Z_k) :=\{ \mu \in P(Z_k) : \mu(\alch_p) =0 \;\forall p
\in \yk(0,0)\}$ and $P_0^+(Z_k) := P_0(Z_k) \cap P^+(Z_k)$.  These are
the weights that are supported on $X_1 \cup X_2$.
It is clear that $\Lk \in P_0(Z_k)$ and that 
$\L \in \W^+\xp \Leftrightarrow \Lk \in P_0^+(Z_k)$.

Similarly given $\g \in \R_s\xp$ define
$$\gk:= \sum_{p \in X_{12} \subset Z_k} \g(p) \, \al_p + \;\;s \!\!\sum_{p
  \in \yk(0,0)} \al_p$$
It is easily seen that (i) $\gk \in Q(Z_k)$ (ii) $\gk \in Q^+(Z_k)
  \Leftrightarrow \g \in R^+_s\xp$ (iii) $\gk \in P_0(Z_k)$ since
  $\g(\xi_1) = \g(\xi_2) =s$.

We now consider two important representation theoretic notions:

Given $\L \in \W^+\xp$, it defines  an integrable highest weight
 representation $L(\Lk)$ of $\kma(Z_k)$ for {\em each} $k \geq 1$. 
Assume $\mu, \nu \in \W^+\xp$.
We can consider the tensor product $L(\Lk) \otimes L(\mk)$; this is
 an integrable representation in category $\mathcal{O}$ and hence
decomposes into a direct sum of integrable highest weight  
representations. This is usually an infinite direct sum, but each
 direct summand occurs with finite multiplicity.
We let $\clmn$ be the multiplicity of the representation
$L(\nu\rk)$ in the decomposition of $L(\Lk) \otimes L(\mk)$.

Similarly given $\beta \in \W\xp$, we consider the 
   {\em branching multiplicity} of the weight $\beta\rk$ for the  action of the
$\kma(A_k)$ corresponding to the subdiagram $A_k$  in the middle. This
number, denoted $\blb(k)$ is defined to be the dimension of the space $\{v \in
L(\Lk)_{\beta\rk}: \mathfrak{n}^+(A_k)\,v=0\}$ i.e vectors of weight $\beta\rk$ 
annihilated
by the positive root spaces of $\kma(A_k)$. 
 
Both the tensor product and branching multiplicities are functions of $k$. If $f: \integers^{>0} \to \integers$ is a function, we say that $f$ {\em stabilizes } if there exists $K$ such that $f(k) = f(k\pr)$ for all $k, k\pr \geq K$. In this case, we set $f(\infty):=f(K)$. We will be interested in conditions under which the $\clmn$ and $\blb(k)$ stabilize. 
\begin{definition}
Given $\L,\mu \in \W\xp$, let $\L \sim \mu$ if $\Lk - \mu\rk \in Q(Z_k)$ for infinitely many values of $k$.
\end{definition}
\begin{remark}\label{remontilda}
It is clear that $\sim$ is an equivalence relation. Further if $(\L +
\mu) \nsim \nu$, then $\clmn$ clearly stabilizes with
$\clmn[\infty]=0$. 
Similarly $\L \nsim \beta$ implies $\blb(k)$ stabilizes with $\blb(\infty)=0$.
\end{remark}

\subsection{Example: the pair $(X(m), A_n)$}\label{oneptfive}
Let $(X, \xi)$ be a given marked \dyn. Fix $m \geq 0$ and set $X_1 =
X(m)$ with $\xi(m)$ distinguished.
Fix $n \geq 0$ and take $X_2 = A_n$ with
the end node being distinguished. Then $Z_k = X(m+n+k)$. This
is the configuration considered in \kv. Weights of the form $\Lk$ for 
$\L \in \W^+\xp$ were called {\em double headed weights} there. 
Under a further hypothesis on $X$ (the so called  
{\em extensibility condition} \cite[definition 2.4]{KV}, see also section \ref{kvreview}
 below), \kv obtained a stabilization 
result for  tensor product multiplicities $\clmn$. We state
this and a related result  in our present notation.

\begin{proposition}\label{kvrestrstat}
Let $X$ be an extensible marked \dyn. For $m,n \geq 0$, let
$X_1:=X(m)$, $X_2:=A_n$, and $Z_k:=Z_k(X_1,X_2)$.
Take $\eta, \delta \in \W\xp$ such that $\eta \sim \delta$.
Then there exists $s \in \integers$ and $\gamma \in \R_s\xp$ such that 
\begin{equation}
\eta\rk - \delta\rk = \gamma\rk \;\; \forall k \geq 1
\end{equation}
\end{proposition}

This is essentially proposition 4.1 of \kv. 
It is here that weights of the form $\gamma\rk$, $\gamma \in \R_s\xp$
appear naturally.

\begin{proposition}\label{kvrestrstat2}
Let $X$ be an extensible marked \dyn. For $m,n \geq 0$, let
$X_1:=X(m)$, $X_2:=A_n$, and $Z_k:=Z_k(X_1,X_2)$.
Suppose $\L, \mu, \nu \in \W^+\xp$ are such 
that there exists $s \in \integers$ and $\gamma \in \R_s\xp$ such that 
\begin{equation}\label{lmng}
\Lk + \mu\rk - \nu\rk = \gamma\rk \;\; \forall k \geq 1
\end{equation}
Then for all $k, k\pr > 2s$,    $\clmn = \clmn[k\pr]$. Thus $\clmn$ stabilizes.
\end{proposition}

Observe that if $\gamma \not\in R^+_s\xp$, then $\gamma\rk \not\in Q^+(Z_k)$ for all $k
\geq 1$. This implies that $\clmn =0 \; \forall k \geq 1$. So we may
as well assume $s \geq 0$ and $\gamma \in R^+_s\xp$.
The proof that $\clmn$ stabilizes,  crucially depends on the fact
that   $\Lk + \mu\rk - \nu\rk$  has this specific form with $\gamma
\in R^+_s\xp$.  This proof appears as part of the proof of theorem 4.5 of \kv in
sections 4.3-4.4 there. 

\subsection{Main theorem}
We now extract the crux of the argument in \kv that proves
Proposition~\ref{kvrestrstat} and formulate a more general
theorem concerning tensor product as well as branching multiplicities.
This works for arbitrary pairs of marked \dyns, but we impose a
condition on our weights that is analogous to equation~\eqref{lmng}.

\begin{theorem}\label{mainthm1}
Let $X_1, X_2$ be arbitrary marked \dyns and $\L, \mu, \nu \in \W^+\xp$,
$\beta \in \W\xp$. 
\be
\item 
Suppose for some $s \geq 0$, $\exists\, \g \in R^+_s\xp$
  with $\Lk + \mu\rk - \nu\rk = \gk$ for all $k \geq 1$. Then for all $k,
  k\pr > 2s$, $\clmn = \clmn[k\pr]$.
\item  Suppose for some $s \geq 0$,  $\exists\, \g \in R^+_s\xp$
  such that $\Lk - \beta\rk = \gk$ for all $k \geq 1$, then for all $k,
  k\pr > 2s$, $\blb(k) = \blb(k\pr)$.
\ee
\end{theorem}

As remarked above, the essential ideas of this proof are the similar
to those of proposition~\ref{kvrestrstat2}; 
We summarize the main steps below (mostly referring to \kv for the
proofs). We also deduce the statement regarding the $\blb$ which does not appear in \kv.

\subsection{Littelmann path model}
As a first step we use the explicit combinatorial description of
$\blb(k)$ and $\clmn$ given by Littelmann's path model. 
For  $k \geq 1$, we let $\Pk$ denote the set of piecewise linear
paths $\pi: [0,1] \to \csa^*(Z_k)$ such that $\pi(0)=0$. To each node
$p \in \nk$, we associate  raising and lowering operators
$e_p, f_p$ on $\integers \Pk$ defined as follows: let 
$\pi \in \Pk$ and $\pi_p(t): = \pi(t)(\Check{\alpha}_p^{(n)})$ for $0 \leq t \leq 1$.
 We consider the function $a: [0,1] \rightarrow [0,1]$ defined by
$a(t) = \min\{1, \pi_p(s) - m_p | t \leq s \leq 1\}$, where $m_p =
 \min\{\pi_p(t)|0 \leq t \leq 1\}$.  Note that $a$ is an increasing
 function.
If $a(1) < 1$ , $f_p \pi
 :=0$. Otherwise, $f_p \pi$ is the path defined by
\begin{equation}\label{fp}
f_p \pi (t) := \pi (t) - a(t) \alpha_p
\end{equation}

Similarly we consider the increasing function 
$b: [0,1] \rightarrow [0,1]$ with
$b(t) = \max\{0, 1 - (\pi_p(s) - m_p) | 0 \leq s \leq t\}$. If $b(0) > 0$,
we set $e_p \pi =0$; otherwise 
\begin{equation}\label{ei}
e_p \pi(t) := \pi (t) + b(t) \alpha_p
\end{equation}
See Littelmann's papers \cite{L1}, \cite{L2}, \cite{L3} for a 
more pictorial desciption of these operators, and
Stembridge\cite{stem} for an axiomatic formulation.

Let $\plk(t):=t\Lk$ be the straight line path with $\plk(1) = \Lk$.
Paths that are obtained by repeated action of the $f_p$, $p \in
Z_k$ on $\plk$ are called Lakshmibai-Seshadri (L-S) paths of shape $\Lk$. 
Consider the following sets:
\begin{align*}
\pk &:= \{\text{ L-S paths } \pi \text{ of shape } \Lk \text{ with }
\pi(1) = \beta\rk \} \\
\pkp &:= \{ \pi \in \pk: e_p \pi =0 \; \forall p \in \yk(0,0)\} \\
\pkm &:= \{\pi \in \pk: \pi \text{ is } \mu\rk \text{ dominant }\} 
\end{align*}
In the last equation, $\mu\rk$ dominance of $\pi$ means that the
shifted path $\mu\rk + \pi$ lies completely in the dominant Weyl chamber.
We now have
\begin{theorem}(Littelmann \cite{L2})
\be
\item $\blb(k) = \# \pkp$
\item $\clmn = \# \pkm[\nu - \mu]$
\ee
\end{theorem}

\subsection{Proof of theorem~\ref{mainthm1}}
We define an auxiliary set. For $k > 2s$ let
$$\Sk := \{ \pi \in \Pk: \pi(t)(\alch_p) = 0 \; \forall t \in [0,1]\, \forall p
\in \yk(s,s)\}$$
These are all (not just L-S) paths that are {\em supported} on the complement of
$\yk(s,s)$.
For all $k, k\pr > 2s$, we have bijections $\phi_{kk\pr}: \Sk \to
\Sk[k\pr]$; given $\pi \in \Sk$, there exist functions $f_p(t)$ such that
\begin{equation}\label{fpeqn}
\pi(t) = \sum_{p \in X_1(s) \cup X_2(s)} f_p(t) \, \om_p
\end{equation}
We define $\phi_{kk\pr}(\pi)$ by the same formula as 
on the right hand side except that we now
interpret $X_1(s)$ and $X_2(s)$ as subdiagrams of $Z_{k\pr}$. It is
clear that $\phi_{kk\pr}$ and $\phi_{k\pr k}$ are inverses of each
other.

The following is the important proposition which details the
relationship between these sets.

\begin{proposition}
Let $k, k\pr > 2s$. Then 
\begin{equation}\label{starthis}
\pkm \subset \pkp \subset \Sk 
\end{equation}
Further, $\phi_{kk\pr}: \Sk \to \Sk[k\pr]$ preserves these
  subsets i.e,
\begin{align}
\phi_{kk\pr}(\pkp) &= \pkppr \label{twothis}\\
\phi_{kk\pr}(\pkm) &= \pkmpr \label{threethis}
\end{align}
\end{proposition}
Proof: The first inclusion in \eqref{starthis} follows
directly from the definition of $\mu\rk$ dominance and the fact that
$\mu\rk(\alch_p)=0$ for all $p \in \yk(0,0)$. 
The proof of the inclusion $\pkp \subset \Sk$ requires
a careful argument with Littelmann paths; this appears in Section 4.4 of
\kv. This latter argument also proves the following fact: 
\begin{equation}\label{lsness}
\pi \in \pkp \Rightarrow \phi_{kk\pr}(\pi) \text{ is an L-S path of
shape } \L^{(k\pr)}
\end{equation}

Given this fact \eqref{lsness} it is easy to  prove \eqref{twothis}
and \eqref{threethis}: we let $\pi \in \pko$ with
$\pi(t) = \displaystyle\sum_{p \in  X_1(s) \cup X_2(s)} f_p(t) \,\om_p $ as in
 equation~\eqref{fpeqn}. We also write $\mu\rk = \sum_{p \in
   X_{12}} b_p \,\om_p$. Then
\begin{align}
\pi \in \pkp &\Leftrightarrow f_p(t) \geq 0 \, \forall t \in [0,1] \, \forall p
\in \yk(0,0) \\
\pi \in \pkm &\Leftrightarrow \pi \in \pkp \text{ and } b_p + f_p(t) \geq 0 \, \forall t
\in [0,1]
\, \forall p \in  X_{12}
\end{align}
These very same conditions ensure that $\pi \in \pkppr$ or $\pi \in
\pkmpr$ as the case may be. \qed

\begin{corollary}
From equation~\eqref{twothis}, we have $\blb(k) =
\blb(k\pr)$. Taking $\beta = \nu - \mu$,
equation~\eqref{threethis} shows $\clmn = \clmn[k\pr]$ for all $k,
k\pr > 2s$.
This  proves our main theorem ~\ref{mainthm1}. \qed
\end{corollary}

\section{Extensible pairs}\label{extpair}
We briefly revisit the situation considered in section
\ref{oneptfive}. Let $X$ be an extensible marked \dyn, $m,n \geq 0$
and set $X_1 :=X(m), X_2 = A_n$ so that $Z_k(X_1,X_2)=X(m+n+k)$. Given
$\L, \mu, \nu \in \W^+\xp$, we have:

\noindent
{\em Case 1:} $(\L + \mu) \nsim \nu$

Remark \ref{remontilda} implies that $\clmn$ stabilizes with $\clmn[\infty]=0$.

\noindent
{\em Case 2:} $(\L + \mu) \sim \nu$

Proposition \ref{kvrestrstat} $\Rightarrow \exists\,\, \gamma \in \R_s\xp$
such that $\Lk + \mu\rk - \nu\rk = \gamma\rk$. Proposition
\ref{kvrestrstat2} then implies stabilization of $\clmn$.
Thus $\clmn$ stabilizes for {\bf all} triples $\L, \mu, \nu \in \W^+\xp$.

For arbitrary $\xp$, Remark \ref{remontilda} still applies. The conclusion of
Proposition \ref{kvrestrstat} however fails in general (see
equation~\eqref{eqabove} for the $B_n$).
So, we can deduce
stabilization of $\clmn$ and $\blb(k)$ only when $\L, \mu, \nu, \beta$
have the special form of Theorem \ref{mainthm1}.

Our next goal will be to define a class of pairs $\xp$
for which Proposition \ref{kvrestrstat} holds. These will be called
{\em extensible pairs}. They include the earlier situation as a
special case; when $X$ is an extensible diagram, $(X(m), A_n)$ will
turn out to be extensible pairs for $m,n \geq 0$. When $\xp$ is an
\ep, one again gets stabilization of $\clmn$ and $\blb(k)$ 
for all $\L, \mu, \nu \in \W^+\xp$, $\beta \in W\xp$.

Sections \ref{extpair} - \ref{stabrepring} below will be concerned exclusively with 
extensible pairs.  We study their properties, define 
their associated {\em number of boxes} functions and finally use the
stable values of $\clmn$ to define a notion of a {\em stable
  representation ring}. 

Among the classical types $A,B,C,D$, only type $A$ falls within the class of extensible pairs.
Readers interested in the $B,C,D$ types may skip directly to section \ref{BCD}.

\subsection{Review of results from \kv}\label{kvreview}
First, we review some relevant notions from \kv.
Given a marked Dynkin diagram  $(X,\xi)$ with $d$ nodes, by a {\em numbering} of $X$
we will mean a bijection $\epsilon : \Nx[X] \to \{1,2,\cdots,d\}$ such
that $\epsilon(\xi) = d$.
For the classical Dynkin diagram $A_k$, $k \geq 1$ we fix the 
numbering $\epsilon_A : \Nx[A_k]
 \to \{1,2,\cdots, k\}$ shown below:
$$
\setlength{\unitlength}{15pt}
\begin{picture}(7,1)(1,-.5)
\put(0,0){\circle*{.25}}
\put(1,0){\circle*{.25}}
\put(2,0){\circle*{.25}}
\put(5,0){\circle*{.25}}
\put(0,0){\line(1,0){3}}
\put(5,0){\line(-1,0){1}}
\put(3.05,0){\makebox(1,0){...}}
\put(0,-.5){\makebox(0,0){\tiny 1}}
\put(1,-.5){\makebox(0,0){\tiny 2}}
\put(2,-.5){\makebox(0,0){\tiny 3}}
\put(5,-.5){\makebox(0,0){\tiny $k$}}
\end{picture}$$
Given a numbering $\epsilon$ of $(X,\xi)$,  the diagram $X(m)$
 inherits a natural numbering $j$
defined by $j(p) = \epsilon(p)$ for $p \in X \subset X(m)$ and
$j(p) = \epsilon_A(p) + d$ for $p \in A_m \subset X(m)$.
\setlength{\unitlength}{15pt}
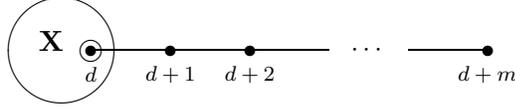
$$\begin{picture}(4,2)(1,-1)
\put(-1,0.25){\makebox(0,0){\bf \large X}}
\put(-.75,0){\circle{3}}
\put(0,0){\circle*{.25}}
\put(0,0){\circle{.5}}
\put(2,0){\circle*{.25}}
\put(4,0){\circle*{.25}}
\put(10,0){\circle*{.25}}
\put(0,0){\line(1,0){6}}
\put(10,0){\line(-1,0){2}}
\put(6.5,0){\makebox(1,0){$\cdots$}}
\put(0,-.6){\makebox(0,0){\footnotesize $d$}}
\put(2,-.6){\makebox(0,0){\footnotesize $d+1$}}
\put(4,-.6){\makebox(0,0){\footnotesize  $d+2$}}
\put(10,-.6){\makebox(0,0){\footnotesize  $d+m$}}
\end{picture}$$
These numberings of $X(m)$ are compatible for different $m$'s, i.e,
 for $m\pr < m$, the numbering $j$ of $X(m)$ restricted to $X(m\pr)
 \subset X(m)$ gives the numbering $j$ of $X(m\pr)$.

Recall that for any \dyn $X$ for which $\det X \neq 0$,  
$P(X)/Q(X)$ is a finite abelian group of order $|\det X|$.
For any $\eta \in P(X)$, we let $[\eta]$ denote its image in
$P(X)/Q(X)$.

Given a marked \dyn $X$
define $$\D[X]:=\det X - \det X(-1)$$
where $X(-1)$ is the diagram obtained from $X$ by deleting the node
$\xi$ and all edges incident on it.
 We then have the
following formula from Equation~(2.2) of \kv:
\begin{equation}\label{detf}
\det X(m) = \det X + m \D[X] \;\; \forall m \geq -1
\end{equation}
If $X=A_1$, then $X(-1)$ is empty; in this case we set $\det X(-1) :=1$. 
\begin{definition}
The marked \dyn $X$ is {\em extensible} if $\det X \neq 0, \D[X] \neq
0$ and $\gcd(\det X, \D[X])=1$.
\end{definition}
It was shown in \kv that the extensibility of $X$ has many pleasant
consequences. Notably, from Lemma~(3.1) of \kv:
\bprop\label{PQcyclic}
Let $X$ be extensible and suppose $m \geq 0$ is such that $\det X(m)
\neq 0$. Then $P(X(m))/Q(X(m))$ is a cyclic group.
Further if we let $\bar{\om}^{(m)}$ denote the fundamental weight corresponding to the
end node of $X(m)$ (i.e $p$ s.t $j(p) = d+m$), then $[\bar{\om}^{(m)}]$
generates this cyclic group.
\eprop
More importantly, Proposition~(3.3) of \kv gives
\begin{proposition}\label{propA}
Let $X$ be extensible.
There exists a unique sequence $(a_i)_{i=1}^\infty$ of integers such that 
 $\forall m \geq 0$ with $\det X(m) \neq 0$ and $\forall p \in
 X(m)$, the relation
\beq\label{starr}
-\D[X] \om_p \equiv a_{j(p)} \, \bar{\om}^{(m)} \pmod{Q(X(m))} 
\eeq
holds in $P(X(m))$.
\end{proposition}
Recall that $j$ is the node numbering introduced above.
It is clear by proposition~\ref{PQcyclic} that $[-\D[X] \om_p ]$
must be a multiple of $[\bar{\om}^{(m)}]$; the content of
proposition~\ref{propA} is that a single sequence $(a_i)$ makes
equation~\eqref{starr} hold {\em for all} values of $m$ under consideration.

We will also need the following lemma which is essentially
Equation~(4.2) of \kv.
\begin{lemma}\label{lemB}
For fixed $m \geq 0$, and $p \in X(m)$, if we write
 $$ -\D[X] \om_p - a_{j(p)} \, \bar{\om}^{(m)} = \sum_{q \in X(m)} 
c_{j(q)} \, \al_q$$ 
with $c_i \in \integers$, then $c_{d+m} = -a_{j(p)}$.
\end{lemma}

\subsection{Definition of extensible pairs}
The goal of this subsection is to generalize the notion of
extensibility to pairs of marked Dynkin diagrams.

Let $(X_1, X_2)$ be a given pair of marked Dynkin diagrams. 
We let $\D[i] :=
\D[X_i] = \det X_i - \det X_i(-1)$. We have
\begin{lemma}\label{detformula}
For all $k \geq 1$
\beq\label{theform}
\det Z_k = (k-1) \D[1] \D[2] + (\det X_1 \, \D[2] + \det X_2  \,
\D[1])
\eeq
\end{lemma}
\noindent
{\em Proof:} We first prove this for $k=1$.
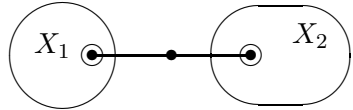
\setlength{\unitlength}{15pt}
$$\begin{picture}(5,3)(1,-1)
\put(-1,0.25){\makebox(0,0){\bf \large $X_1$}}
\put(-.75,0){\circle{3.2}}
\put(0,0){\circle*{.25}}
\put(0,0){\circle{.5}}
\put(2,0){\circle*{.25}}
\put(4,0){\circle*{.25}}
\put(4,0){\circle{.5}}
\put(4.75,0){\oval(3.5,2.5)}
\put(5.5,0.5){\makebox(0,0){\bf \large $X_2$}}
\put(0,0){\line(1,0){4}}
\end{picture}$$
Put $\bar{X}_i := X_i(-1)$ for $\ri$. Now,
the matrix $C(Z_1)$ looks like
\begin{equation*}
C(Z_1) = \left[
\begin{array}{r|c|l}
C(X_1)&&\\
      &-1&\\
\hline
    -1&2&-1\\
\hline
      &-1&\\
     &&C(X_2)\\
\end{array}
\right]
\end{equation*}
Expanding along the middle column we get
\begin{equation}\label{detz1}
\det Z_1  = 2 \det X_1 \det X_2  + \det P_1  + \det P_2
\end{equation}
where $P_1$ is a matrix of the form
\begin{equation*}
P_1 = \left[
\begin{array}{r|c|l}
C(\bar{X}_1)&*&\\
      &*&\\
\hline
    &-1&-1\\
\hline
&&\\
    &&C(X_2)
\end{array}
\right]
\end{equation*}
The upper triangularity of $P_1$ gives 
$\det P_1  = - \det \bar{X}_1  \det X_2 $. A similar calculation shows
that $P_2$ is a lower triangular matrix with
$\det P_2  = - \det \bar{X}_2  \det X_1 $. Putting these back in
Equation\eqref{detz1} completes the proof for $k=1$. 

If $k>1$, then consider $X_1(k-1)$.
We can think of $Z_k = Z_k(X_1, X_2)$ as being the same thing as
$Z_1(X_1(k-1), X_2)$. In terms of a picture:
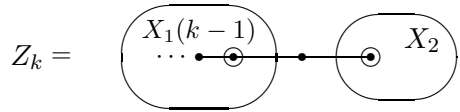
\setlength{\unitlength}{13pt}
$$\begin{picture}(7,3)(-2,-1)
\put(-6,0){\makebox(0,0){\bf \large $Z_k$}}
\put(-5,0){\makebox(0,0){=}}

\put(-1,0.75){\makebox(0,0){\bf $X_1(k-1)$}}
\put(-1,0){\oval(4.5,3)}
\put(0,0){\circle*{.25}}
\put(0,0){\circle{.5}}
\put(-1,0){\circle*{.25}}
\put(2,0){\circle*{.25}}
\put(4,0){\circle*{.25}}
\put(4,0){\circle{.5}}
\put(4.75,0){\oval(3.5,2.5)}
\put(5.5,0.5){\makebox(0,0){\bf \large $X_2$}}

\put(0,0){\line(1,0){4}}
\put(0,0){\line(-1,0){1}}
\put(-1.75,0){\makebox(0,0){$\cdots$}}
\end{picture}$$
Applying our $k=1$
result gives 
$$\det Z_k  = \det X_1(k-1) \, \D[2]+ \det X_2 \,(\det X_1(k-1) -\det X_1(k-2))$$
Lemma~\ref{detformula} implies
$\det X_i(m) = \det X_i + m\D[i] \;\; \forall m \geq -1$.
Putting this back we get:
$$ \det Z_k = (\det X_1 + (k-1)\D[1]) \,\D[2]+ \det X_2 \,\D[1]$$
proving our lemma. \qed

\begin{definition}\label{ep}
A pair of marked \dyns $(X_1, X_2)$ is called  an {\em extensible pair} if
\be
\item Each $X_i$ is an extensible diagram, i.e, $\gcd(\det X_i, \D[i])
  =1$, $\D[i] \neq 0$, $\det X_i \neq 0$, $i=1,2$.
\item $\gcd(\D[1],  \D[2])=1$.
\ee
\end{definition}

\begin{remark}
\be
\item Observe that $A_1$ is an extensible \dyn since $\D[A_1]=1$. It
  is clear that $X$ is extensible $\iff$ $(X, A_1)$ is an extensible
  pair.
\item If $k_0$ is such that $\det Z_{k_0} =0$, then many of our
  familiar statements for $\kma(Z_{k_0})$ break down; 
for instance $P(Z_{k_0})/Q(Z_{k_0})$ is no  longer a finite group, 
the $\om_p$'s do not span $\csa^*(Z_{k_0})$ etc.
But if \xp is an extensible pair, $\D[1] \D[2] \neq 0$. By
  Lemma~\ref{detformula} there can be at most one $k_0$ for
  which $\det Z_{k_0} =0$. So this exceptional situation can occur
  for at most one value of $k$. Most of our later results will only 
hold  for $k \neq k_0$.
\ee
\end{remark}

\begin{lemma}\label{newlem}
If \xp is an extensible pair, then $\gcd(\D[1]\D[2], \det Z_k) =1 \;
\forall k \geq 1$
\end{lemma}
\noindent
{\bf Proof:} Let $P$ be a prime such that $P|\,\D[1]\D[2]$. Then $P |
\,\D[1]$ or $P|\,\D[2]$. Suppose $P|\,\D[1]$, definition~\ref{ep} implies
that $P$ does not divide either $\det X_1$ or $\D[2]$. 
From Equation~\eqref{theform}, this means
that $P$ cannot divide $\det Z_k$. The $P | \, \D[2]$ case is similar.\qed

\begin{remark}\label{exceptions}
It is easily seen that if  $\xp$ is an \ep, then so is $(X_1(m),
X_2(n))$ for almost all $m,n \geq 0$. The only exceptional values are
those which make $\det(X_1(m))=0$ or $\det(X_2(n))=0$.
\end{remark}

\begin{example}
Let $\mathcal{U}$ be the set of marked Dynkin diagrams $U$ with $\det
U \neq 0$ and $\Delta_U = \pm 1$. From  Table 1 of \kv,
it is clear that Types $A, E, F^{(1)},  F^{(2)},  G^{(1)},  G^{(2)}$ are
all in $\mathcal{U}$. 

Given $X_1, X_2 \in \mathcal{U}$, it is clear that $(X_1, X_2)$ is an
extensible pair. Some of the sequences $Z_k$ obtained thus are :
\begin{enumerate}
\renewcommand{\theenumi}{\roman{enumi}}
\item $X_1 = X_2 = A_1$. 
\setlength{\unitlength}{15pt}
\begin{picture}(5,.5)(-5,0)
\put(-1,0){\makebox(0,0){$Z_k$ = }}
\put(1,0){\circle*{.25}}
\put(2,0){\circle*{.25}}
\put(3,0){\circle*{.25}}
\put(4,0){\circle*{.25}}
\put(7.5,0){\circle*{.25}}
\put(8.5,0){\circle*{.25}}
\put(9.5,0){\circle*{.25}}
\put(1,0){\line(1,0){3.5}}
\put(9.5,0){\line(-1,0){2.5}}
\put(5.3,0){\makebox(1,0){$\cdots$}}

\end{picture}
\vspace{0.6cm}
\item $X_1 = E_6, X_2 = A_1$
\setlength{\unitlength}{15pt}
\begin{picture}(5,1)(-5,0)
\put(-2,0){\makebox(0,0){$Z_k$ = }}
\put(0,0){\circle*{.25}}
\put(1,0){\circle*{.25}}
\put(2,0){\circle*{.25}}
\put(2,1){\circle*{.25}}

\put(3,0){\circle*{.25}}
\put(4,0){\circle*{.25}}
\put(7.5,0){\circle*{.25}}
\put(8.5,0){\circle*{.25}}
\put(9.5,0){\circle*{.25}}
\put(0,0){\line(1,0){4.5}}
\put(2,0){\line(0,1){1}}
\put(9.5,0){\line(-1,0){2.5}}
\put(5.3,0){\makebox(1,0){$\cdots$}}

\end{picture}
\vspace{0.6cm}
\item $X_1 = X_2 = E_6$
\setlength{\unitlength}{15pt}
\begin{picture}(5,1)(-6.2,0)
\put(-2,0){\makebox(0,0){$Z_k$ = }}
\put(0,0){\circle*{.25}}
\put(1,0){\circle*{.25}}
\put(2,0){\circle*{.25}}
\put(2,1){\circle*{.25}}

\put(3,0){\circle*{.25}}
\put(4,0){\circle*{.25}}
\put(7.5,0){\circle*{.25}}
\put(8.5,0){\circle*{.25}}
\put(8.5,1){\circle*{.25}}
\put(9.5,0){\circle*{.25}}
\put(10.5,0){\circle*{.25}}
\put(0,0){\line(1,0){4.5}}
\put(2,0){\line(0,1){1}}
\put(10.5,0){\line(-1,0){3.5}}
\put(8.5,0){\line(0,1){1}}
\put(5.3,0){\makebox(1,0){$\cdots$}}

\end{picture}
\vspace{0.6cm}
\item $X_1 = E_6, X_2 = G_2$
\setlength{\unitlength}{15pt}
\begin{picture}(5,1)(-5,0)
\put(-2,0){\makebox(0,0){$Z_k$ = }}
\put(0,0){\circle*{.25}}
\put(1,0){\circle*{.25}}
\put(2,0){\circle*{.25}}
\put(2,1){\circle*{.25}}

\put(3,0){\circle*{.25}}
\put(4,0){\circle*{.25}}
\put(7.5,0){\circle*{.25}}
\put(8.5,0){\circle*{.35}}
\put(9.5,0){\circle*{.35}}
\put(0,0){\line(1,0){4.5}}
\put(2,0){\line(0,1){1}}
\put(9.5,0){\line(-1,0){2.5}}
\put(9.5,0.1){\line(-1,0){1}}
\put(9.5,-0.1){\line(-1,0){1}}

\put(5.3,0){\makebox(1,0){$\cdots$}}
\put(9,0){\makebox(0,0){\small $>$}}

\end{picture}

\end{enumerate}

\end{example}

\subsection{Node numbering}
Let $(X_1, X_2)$ be an extensible pair and assume $X_i$ has $d_i$
nodes $i=1,2$. We assume we are given
numberings $\epsilon_i$ of $X_i$ ($i=1,2$) as in section \ref{kvreview}.
We use the
numberings $\epsilon_1, \epsilon_2$ and $\epsilon_A$ to construct 
two numberings $i, \bi$ of $Z_k$. Define 
$i, \bi : N(Z_k) \to \{1 \leq m \leq d_1 + d_2 + k\}$ by
\begin{equation*}
i(p):=
\begin{cases}
\epsilon_1(p) &\text{if $p \in N(X_1)$}\\
\epsilon_A(p) + d_1 &\text{if $p \in N(A_k)$}\\
(d_1 + d_2 + k+ 1)-\epsilon_2(p) &\text{if $p \in N(X_2)$}
\end{cases}
\end{equation*}
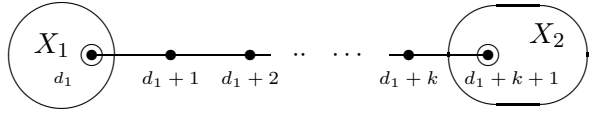
\setlength{\unitlength}{15pt}
$$\begin{picture}(5,3)(2,-1)
\put(-1,0.25){\makebox(0,0){\bf \large $X_1$}}
\put(-.75,0){\circle{3.2}}
\put(0,0){\circle*{.25}}
\put(0,0){\circle{.5}}
\put(2,0){\circle*{.25}}
\put(4,0){\circle*{.25}}
\put(8,0){\circle*{.25}}
\put(10,0){\circle*{.25}}
\put(10,0){\circle{.5}}
\put(10.75,0){\oval(3.5,2.5)}
\put(11.5,0.5){\makebox(0,0){\bf \large $X_2$}}
\put(0,0){\line(1,0){4.5}}
\put(10,0){\line(-1,0){2.5}}
\put(5,0){\makebox(0.5,0){$\cdot\cdot$}}
\put(6,0){\makebox(1,0){$\cdots$}}
\put(-.7,-.6){\makebox(0,0){\tiny $d_1 $}}
\put(2,-.6){\makebox(0,0){\scriptsize $d_1 + 1$}}
\put(4,-.6){\makebox(0,0){\scriptsize  $d_1+2$}}
\put(8,-.6){\makebox(0,0){\scriptsize  $d_1+k$}}
\put(10.6,-.6){\makebox(0,0){\scriptsize  $d_1+k+1$}}
\end{picture}$$

\begin{equation*}
\bi(p):=
\begin{cases}
\epsilon_2(p) &\text{if $p \in N(X_2)$}\\
(d_2 + k + 1) - \epsilon_A(p) &\text{if $p \in N(A_k)$}\\
(d_1 + d_2 + k+ 1)-\epsilon_1(p) &\text{if $p \in N(X_1)$}
\end{cases}
\end{equation*}
Our earlier figure showed $Z_k$ with the numbering given by
$i(\cdot)$. The figure below shows the numbering via $\bi(\cdot)$.
\setlength{\unitlength}{15pt}
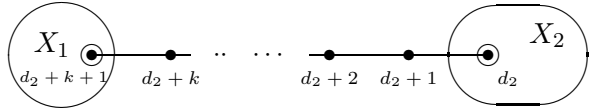
$$\begin{picture}(5,3)(2,-1)
\put(-1,0.25){\makebox(0,0){\bf \large $X_1$}}
\put(-.75,0){\circle{3.2}}
\put(0,0){\circle*{.25}}
\put(0,0){\circle{.5}}
\put(2,0){\circle*{.25}}
\put(6,0){\circle*{.25}}
\put(8,0){\circle*{.25}}
\put(10,0){\circle*{.25}}
\put(10,0){\circle{.5}}
\put(10.75,0){\oval(3.5,2.5)}
\put(11.5,0.5){\makebox(0,0){\bf \large $X_2$}}
\put(0,0){\line(1,0){2.5}}
\put(10,0){\line(-1,0){4.5}}
\put(3,0){\makebox(0.5,0){$\cdot\cdot$}}
\put(4,0){\makebox(1,0){$\cdots$}}
\put(-.7,-.6){\makebox(0,0){\tiny $d_2 + k +1$}}
\put(2,-.6){\makebox(0,0){\scriptsize $d_2 + k$}}
\put(6,-.6){\makebox(0,0){\scriptsize  $d_2+2$}}
\put(8,-.6){\makebox(0,0){\scriptsize  $d_2+1$}}
\put(10.5,-.6){\makebox(0,0){\scriptsize  $d_2$}}
\end{picture}$$

\subsection{$P(Z_k)/ Q(Z_k)$ is cyclic}
Let $\bri :=\yk(0,0) \cup \{\xi_1, \xi_2\}$.
Our new notion of {\em extensible pairs} has the same nice consequence as the
previous notion of {\em extensible diagrams}:
\begin{lemma}\label{pqcyclic}
Let \xp be an extensible pair. For $k \geq 1$, if $\det Z_k \neq 0$,
then $P(Z_k)/ Q(Z_k)$ is a cyclic group. Further, if $u,v \in \bri
\subset \nk$ are any two adjacent nodes then $[\om_u - \om_v]$ generates this group.
\end{lemma}
\noindent
{\bf Proof:} Let $p, p^{\prime},  p^{\prime\prime} \in \bri$ be 3
adjacent nodes as in figure 
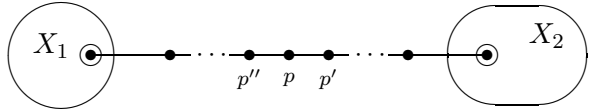
\setlength{\unitlength}{15pt}
$$\begin{picture}(5,3)(3,-1)
\put(-1,0.25){\makebox(0,0){\bf \large $X_1$}}
\put(-.75,0){\circle{3}}
\put(0,0){\circle*{.25}}
\put(0,0){\circle{.5}}
\put(2,0){\circle*{.25}}
\put(4,0){\circle*{.25}}
\put(5,0){\circle*{.25}}
\put(6,0){\circle*{.25}}
\put(8,0){\circle*{.25}}
\put(10,0){\circle*{.25}}
\put(10,0){\circle{.5}}
\put(10.75,0){\oval(3.5,2.5)}
\put(11.5,0.5){\makebox(0,0){\bf \large $X_2$}}
\put(0,0){\line(1,0){2.5}}
\put(10,0){\line(-1,0){2.5}}
\put(3.5,0){\line(1,0){3}}
\put(2.6,0){\makebox(1,0){$\cdots$}}
\put(6.6,0){\makebox(1,0){$\cdots$}}
\put(4,-.6){\makebox(0,0){\footnotesize $p\ppr$}}
\put(5,-.6){\makebox(0,0){\footnotesize $p$}}
\put(6,-.6){\makebox(0,0){\footnotesize  $p\pr$}}
\end{picture}$$

Then $[\om_p - \om_{p\pr}] = [\om_{p\ppr} -
  \om_p]$ in $P(Z_k) / Q(Z_k)$ since
$$ (\om_p - \om_{p\pr}) -  
(\om_{p\ppr} - \om_p) =  2 \om_p - \om_{p\pr} - \om_{p\ppr} = \al_p
\in \qzk$$
Iterating this argument, it is clear that the elements 
$[\om_u - \om_v]$ with $u,v \in \bri$ and $i(v) = i(u) + 1$ are all
equal. Further if $[\om_u - \om_v]$ is a generator of the group, then
so is $[\om_v - \om_u]$. It is thus enough to prove the lemma for a
fixed choice of $u,v \in \bri$ with $i(v) = i(u) + 1$.

We write $\om_u - \om_v  = \sum_{p \in Z_k}  c_p \al_p$ with $c_p
\in \rationals$.  We now compute $c_u$ and $c_v$. Set $a = i(u) - d_1
- 1$ and $b = \bi(v) - d_2 -1$, so $a+b = k-2$.

The coefficient of $\al_u$ in $\om_u$ is the $(u,u)^{th}$ element of
$C(Z_k)^{-1}$, the inverse of the generalized Cartan matrix. This
coefficient is thus:
$$\frac{\text{ cofactor of }(u,u)^{th} \text{ element of }
  C(Z_k)}{\det Z_k}
= \frac{1}{\det Z_k} \det X_1(a) \det X_2(b+1)$$

Similarly the coefficient of $\al_v$ in $\om_u$ is the $(v,u)^{th}$ element of
$C(Z_k)^{-1}$, and is thus $ -\det A / \det Z_k$ where 

\begin{equation*}
A = \left[
\begin{array}{r|c|l}
C(X_1(a))&&\\
      &-1&\\
\hline
    &-1&-1\\
\hline
     &&C(X_2(b))\\
\end{array}
\right]
\end{equation*}
The upper triangularity of $A$ gives $\det A = - \det X_1(a)
\det X_2(b)$. We now switch the  roles of $u$ and $v$;
the coefficient of $\al_u$ in $\om_v$  is $\det X_1(a) \det X_2(b)/
\det Z_k$ and coefficient of $\al_v$ in $\om_v$ is  $\det X_1(a+1)
\det X_2(b)/ \det Z_k$. Finally we also have Equation~\eqref{detf} which gives
$\det X_i(r) = \det X_i + r \D$. Putting these all together, we get
\begin{equation}\label{cucv}
c_u - c_v = \frac{\det X_1 \D[2] + \det X_2 \D[1]
+ (k-2) \D[1] \D[2]}{\det Z_k}
\end{equation}
Lemma~\ref{newlem} implies that the numerator and denominator of
Equation~\eqref{cucv} are relatively prime. So, the least
$n \in \N$ such that $n(c_u - c_v) \in \integers$ is $n =
|\det Z_k|$. In particular this implies that the order of the element
$[\om_u - \om_v]$ in $P(Z_k)/Q(Z_k)$ is $\geq |\det Z_k|$. Since 
$\# P(Z_k)/Q(Z_k) = |\det Z_k|$, this completes the proof of
Lemma~\ref{pqcyclic}. \qed

Now fix $u,v \in \bri \subset Z_k$ with $i(v) = i(u) + 1$. Our
next goal is to explicitly write each $[\om_p]$, $p \in Z_k$
  as a multiple of $[\om_u  - \om_v]$.
 for each $p \in Z_k$. This is analogous to Proposition \ref{propA}. 
In fact we will use this latter proposition to deduce our result below.

Consider the subdiagrams $\tx[1]:=X_1(i(u))$ and
$\tx[2]:=X_2(\bi(v))$ of $Z_k$. The 
numberings $i, \bi$ of $Z_k$ can be restricted to $\tx[i] ,\, \ri$ to
give {\em functions} from $\tx[i]$ into the set $\{1,2,\cdots, d_1
+ d_2 + k\}$. Since $(X_1, X_2)$ is an \ep, each $X_i$ is an
extensible diagram. For fixed $i \in \{1,2\}$, let 
$\{\tom_p : p \in \tx[i]\}$ and $\{\tal_p : p \in \tx[i]\}$ 
denote the fundamental weights and simple roots of
$\tx[i]$. We can now apply Proposition~\ref{propA} above to the
$X_i$, $i=1,2$. This gives:
\begin{proposition}\label{propApr}
There exist infinite sequences $(\aip[]{1})_{i=1}^{\infty}$ and 
$(\aip[]{2})_{i=1}^{\infty}$ (determined uniquely by $X_1$ and $X_2$)
such that
\begin{align}
-\D[1]\, \tom_p &\equiv \aip{1} \, \tom_u \pmod{Q(\tx[1])} \;\;\;\;\forall
p \in \tx[1] \\
-\D[2]\, \tom_p &\equiv \aipb{2} \, \tom_v \pmod{Q(\tx[2])}  \;\;\;\; \forall
p \in \tx[2]
\end{align}
\end{proposition}
Additionally, 
\begin{lemma}\label{lembpr}
\be
\item For $p \in \tx[1]$, let
$$ \tilde{\beta}_p:= -\D[1]\, \tom_p - \aip{1} \; \tom_u
= \sum_{q \in \tx[1]} \tilde{c}_{q,p} \, \tal_q$$ with 
$\tilde{c}_{q,p} \in \integers$. Then $\tilde{c}_{u,p} = -\aip{1}$.
\item For $p \in \tx[2]$, let
 $$ \tilde{\beta}_p:= -\D[2] \,\tom_p - \aipb{2} \; \tom_v
= \sum_{q \in \tx[2]} \tilde{c}_{q,p} \, \tal_q$$ with 
$\tilde{c}_{q,p} \in \integers$. Then $\tilde{c}_{v,p} = -\aipb{2}$.

\ee

\end{lemma}
\noindent
{\em Proof:} Follows from lemma~\ref{lemB} gives. \qed

We now consider the original diagram $Z_k$. Fix $p \in \tx[1]$,
take the coefficients $\tilde{c}_q$ given by Lemma~\ref{lembpr} and
consider the element $\beta_p \in Q(Z_k)$ given by 
$$\beta_p := \sum_{q \in \tx[1]} \tilde{c}_{q,p} \,\al_q$$
 The $\al_q$ are
now the simple roots of $Z_k$. Looking at how $\tx[1]$ sits inside
 $Z_k$ as a subdiagram, it is clear that the following relations hold:
\begin{align}
\tilde{\beta}_p (\Check{\tilde{\alpha}}_r) &= \beta_p (\alch_r)
& &\forall r \in \tx[1] \subset Z_k \notag\\
\beta_p(\alch_v) &= - \tilde{c}_{u,p} = \aip{1} & &\label{bet}\\
\beta_p(\alch_q)&=0 & &\forall q \not\in \tx[1] \cup \{v\} \notag
\end{align}
These follow from the fact that
 $u$ is the only node connected to $v$  and from Lemma~\ref{lembpr}.

Equation~\eqref{bet} and the definition of $\tilde{\beta}_p$ imply the
following expression for $\beta_p$ as a linear combination of the
$\om_r, \; r \in Z_k$
\beq\label{eqone}
\beta_p= -\D[1] \om_p - \aip{1} (\omuv) \;\; \forall p \in Z_k
\text{ with }
i(p) \leq i(u)
\eeq

The corresponding picture for $p \in  \tx[2]$ is obtained
similarly; we define
$\beta\pr_p \in Q(Z_k)$ by $$\beta\pr_p := \sum_{q \in \tx[2]} \tilde{c}_{q,p}
\al_q$$
The earlier argument for $\beta_p$ can be carried out with obvious
modifications and gives the following expression for $\beta\pr_p$:
\beq\label{eqonepr}
\beta\pr_p= -\D[2] \om_p - \aipb{2} (\omvu) \;\; \forall p \in Z_k
\text{ with } \bi(p) \leq \bi(v)
\eeq

We now have the following proposition which follows from
Equations~\eqref{eqone}, \eqref{eqonepr} and the fact that $\beta_p,
\beta\pr_p \in Q(Z_k)$.
\bprop\label{congs}
Let $p \in \nk$. Then
\be
\item If $i(p) \leq i(u)$, then  $-\D[1] \om_p \equiv \aip{1} (\omuv)
  \pmod{Q(Z_k)}$
\item If $\bi(p) \leq \bi(v)$, then  $-\D[2] \om_p \equiv \aipb{2}(\omvu)
  \pmod{Q(Z_k)}$
\ee
\eprop

\section{Stabilization of $\clmn$ and $\blb(k)$}\label{h2plus}
\subsection{Two sided dominant weights}\label{notats}
 We will now index dominant integral weights of $\kma(Z_k)$ as in \kv, 
in a two-sided (or {\em double headed}) fashion. Unlike in
section \ref{twofour}, we consider weights whose supports are not necessarily
contained in $X_1 \cup X_2$.

Let
$$ \ho := \{ (x_1,x_2,\cdots) : x_i 
\in \integers \, \forall i \; and \; x_i \neq 0
 \mbox{ for only finitely many }  i \}$$
and
$$ \hop = \{ (x_1,x_2,\cdots) \in \ho : x_i \geq 0 \; \forall i\}$$ 
Given $x=(x_1,x_2,\cdots) \in \ho$ we 
define the {\em length} of $x$ to be: $\ell(x):= \max\{i: x_i \neq 0\}$. 
Let 
\begin{align*}
\Ht &:=\ho \times \ho & \text{and}& & \Htp&:=\hop \times \hop
\end{align*}
Given elements
$x = (x_1,x_2,\cdots)$, $y = (y_1,y_2,\cdots)$ of $\ho$, let
$\gamma=(x,y) \in \Ht$. We define the {\em left support} of $\gamma$
to be $\ls(\gamma) : = \max(\ell(x), d_1)$,  the {\em right support}
as $rs(\gamma) :=\max(\ell(y), d_2)$ and the {\em support}
$\supp(\gamma):=\ls(\gamma) + \rs(\gamma)$. Recall that $d_i$ here is the
number of nodes in $X_i$.
Now, $\gamma$ can be used to define a weight of $\mathfrak{g}(Z_k)$
for all large $k$; specifically for 
$k \geq \supp(\gamma) - d_1 - d_2$ as follows. Define
$$\gamma^{(k)} := \sum_{p \in Z_k} (\xip + \yipb) \, \om_p$$ 
Note that since $k$ is large, 
at most one of $\xip$ or $\yipb$ can be nonzero for each $p \in \nk$. 
The figure below shows each node $p$ of $Z_k$ labeled by the corresponding
$x_i$ or $y_i$.
\setlength{\unitlength}{15pt}
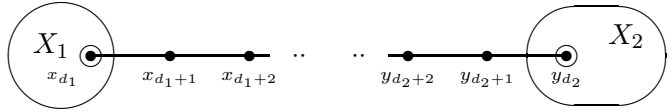
$$\begin{picture}(5,3)(3,-1)
\put(-1,0.25){\makebox(0,0){\bf \large $X_1$}}
\put(-.75,0){\circle{3.2}}
\put(0,0){\circle*{.25}}
\put(0,0){\circle{.5}}
\put(2,0){\circle*{.25}}
\put(4,0){\circle*{.25}}
\put(8,0){\circle*{.25}}
\put(10,0){\circle*{.25}}
\put(12,0){\circle*{.25}}
\put(12,0){\circle{.5}}
\put(12.75,0){\oval(3.5,2.5)}
\put(13.5,0.5){\makebox(0,0){\bf \large $X_2$}}
\put(0,0){\line(1,0){4.5}}
\put(12,0){\line(-1,0){4.5}}
\put(5,0){\makebox(0.5,0){$\cdot\cdot$}}
\put(6.5,0){\makebox(0.5,0){$\cdot\cdot$}}
\put(-.7,-.6){\makebox(0,0){\tiny $x_{d_1}$}}
\put(2,-.6){\makebox(0,0){\scriptsize $x_{d_1+1}$}}
\put(4,-.6){\makebox(0,0){\scriptsize  $x_{d_1+2}$}}
\put(8,-.6){\makebox(0,0){\scriptsize  $y_{d_2+2}$}}
\put(10,-.6){\makebox(0,0){\scriptsize  $y_{d_2+1}$}}
\put(12,-.6){\makebox(0,0){\scriptsize  $y_{d_2}$}}
\end{picture}$$

\subsection{The {\em number of boxes} function}
Given an extensible \dyn $X$, section 3.3 of \kv introduced the so
called {\em number of boxes} function $|\cdot|_X : \Ht \to
\integers$. For an extensible pair \xp, we now define an analogous
function $\nob[\cdot]$. We shall refer to this also as the {\em
  number of boxes} function.
\begin{definition}\label{nobox}
Given $\gamma =(x,y) \in \Ht$ as above, define
$$ \nob := \D[2] \,(\sum_i x_i \aip[]{1}) - \D[1] \,(\sum_i y_i
\aip[]{2})$$
\end{definition}
\begin{remark}
When $X_2 = A_1$ we have $\D[2]=1$ and $\aip[]{2} = i$.
Then $|\gamma|_{(X_1, A_1)}$ coincides with the function
$|\gamma|_{X_1}$ in definition~(3.12) of \kv. 
\end{remark}
To see the significance of $\nob$, we consider the following situation. Let $\gamma
\in \Ht$ and $k \geq \supp(\gamma) - d_1 - d_2$ be fixed. Pick $u, v
\in \bri$  such that $i(u) \geq \ls(\g)$,
 $\bi(v) \geq \rs(\g)$ and $i(v) = i(u) + 1$. Consider $\gk \in
\pzk$. Recall from Lemma~\ref{pqcyclic} that $\gk$
must be congruent modulo $\qzk$ to an integer multiple of $\omuv$.
To calculate this integer explicitly, 
\begin{align}\label{impeqn}
-\D[1]\D[2] \gk - \nob (\omuv) = &\D[2]\left( \sum_{\substack{p \in
    Z_k\\ i(p) \leq i(u)}} \xip \left[ -\D[1]\,\om_p - \aip{1} (\omuv)\right]
\right)\notag\\
                    &+ \D[1]\left( \sum_{\substack{p \in
    Z_k\\ \bi(p) \leq \bi(v)}} \yipb \left[ -\D[2]\,\om_p - \aipb{2} (\omvu)\right]
\right)
\end{align}

Proposition~\ref{congs} shows that both terms within square brackets
in the above equation are elements of $\qzk$. So, we have

\bprop
With notation as above, 
\beq\label{gamcong}
 - \D[1]\D[2] \gk \equiv \nob (\omuv) \pmod{\qzk} 
\eeq
\eprop

Observe that this equation  enables us to write $\gk$ itself as a
multiple of $(\omuv)$ modulo $\qzk$. To see this, note that Lemma~\ref{newlem}
gives
$\gcd(\D[1]\D[2], \det Z_k)=1$. Since $\# \pzk / \qzk = |\det Z_k|$,
we have

\beq\label{kspec}
\gk \equiv  (\D[1]\D[2])^{-1} \nob \;(\omuv) \pmod{\qzk} 
\eeq
where $(\D[1]\D[2])^{-1}$ denotes the inverse of $\D[1]\D[2]$ in
$\frac{\integers}{(\det Z_k)\integers}$.  We also
have the following corollary to the above proposition:
\begin{corollary}\label{corr}
Let $\g \in \Ht$. Then 
\begin{align*}
\nob =0 &\iff \gk \in \qzk \text{ for all large } k \\
&\iff \gk \in \qzk \text{ for infinitely many values of } k
\end{align*}
\end{corollary}

\noindent
{\bf Proof:} Since $\omuv$ generates the cyclic group $\pzk/\qzk$,
Equation~\eqref{kspec} implies that 
$\gk \in \qzk$ {\em iff} $\det Z_k$ divides $\nob$. But as $k \to \infty$,
$|\det Z_k| \to \infty$ since $\D[1]\D[2] \neq 0$; this follows from
Equation~\eqref{theform}. The corollary follows. \qed

\subsection{Depth}

In this subsection, we introduce an important statistic called the
{\em depth} on the set $\{\gamma: \nob=0\}$.

Let $\g=(x,y) \in \Ht$ with $\nob =0$. Fix $k \geq \supp(\gamma) - d_1 -
d_2$ and pick $u, v$ as in previous the subsection.

Consider Equation~\eqref{impeqn}. Imposing the condition $\nob =0$ 
reduces the left hand side to $-\D[1]\D[2] \gk$. As already remarked,
the sum on the right hand side of
\eqref{impeqn} is an element of $\qzk$. We now calculate the
coefficient of $\al_u$ in this sum. Using Lemma~\ref{lembpr} and
Equation~\eqref{eqone}, the coefficient of $\al_u$ on the RHS equals
$$ - \D[2] ( \sum_{\substack{p \in
    Z_k\\ i(p) \leq i(u)}} \xip \aip{1}) = -\D[2] \sum_{i} x_i
    \aip[]{1}$$
Similarly, employing Equation~\eqref{eqonepr}, the coefficient of
    $\al_v$ on the RHS of \eqref{impeqn} becomes
$$ - \D[1] ( \sum_{\substack{p \in
    Z_k\\ \bi(p) \leq \bi(u)}} \yipb \aipb{2}) = -\D[1] \sum_{i} y_i
    \aip[]{2}$$
Observe now that since $\nob =0$, the  coefficients of $\al_u$ and
    $\al_v$ are in fact equal !  Putting everything together we get

\begin{lemma}\label{star}
If $\nob=0$, then we have:
\begin{align}\label{coco}
\text{coefficient of } \al_u \text{ in } & -\D[1]\D[2] \gk 
=  \;\;\text{coefficient of } \al_v \text{ in } -\D[1]\D[2] \gk \notag \\
&=  -\D[2] \sum_{i} x_i \aip[]{1} \; = -\D[1] \sum_{i} y_i \aip[]{2}
\end{align}
\end{lemma}
This motivates the following
\begin{definition}
If $\g = (x,y) \in \Ht$ with $\nob=0$, the {\em depth} of $\g$ is
defined to be
$$ \dep(\g):= \frac{\sum_{i} x_i \aip[]{1}}{\D[1]} =\frac{\sum_{i}
  y_i\aip[]{2}}{\D[2]}$$
\end{definition}
\begin{remark}\label{staremk}
\be
\item Observe that $\dep(\g) \in \integers$ since $\gcd(\D[1], \D[2])=1$.
\item From equation~\eqref{coco}, it is clear that $\dep(\g)$ is
  just the coefficient of $\al_u$ (or $\al_v$) in $\gk$.
\item When $X_2 = A_1$, our notion of depth reduces to the notion
  introduced in \kv, Definition (4.3).
\ee
\end{remark}

Corollary~\ref{corr} guarantees that if $\nob=0$, then $\gk$ can be
written as an integer linear combination of $\{ \al_p: p \in Z_k\}$ for
all large $k$. In this subsection, we study how the coefficients of
this linear combination change with $k$.

In order to be able to refer more easily to nodes of $Z_k$ which occur to the
left/right/middle of the diagram, we define some notation; for $l, r
\in \N$ we let 
\begin{align*}
L_k(l) &:= \{p \in Z_k : i(p) \leq l\} \\
R_k(r) &:= \{p \in Z_k : \bi(p) \leq r\} \\
M_k(l,r) &:= Z_k - (L_k(l) \cup R_k(r))
\end{align*}

We now have the following proposition which generalizes
Proposition~(4.1) of \kv.
\bprop\label{aaa}
Suppose $\gamma \in \Ht$ with $\nob =0$; let $l :=\ls(\g),
r:=\rs(\g)$. Then there exist integers $b_i \,(1 \leq i \leq l-1)$, 
$c_i \,(1 \leq i \leq r-1)$ and $s$ such that for all $k \geq
l+r-d_1-d_2$ we have
\beq\label{bicis}
\gk = \sum_{p \in L_k(l-1)} b_{i(p)} \,\al_p + \sum_{p \in M_k(l-1,
  r-1)} s\,\al_p +  \sum_{p \in R_k(r-1)} c_{\bi(p)} \,\al_p
\eeq
\eprop
\noindent
{\bf Proof:} Let $k$ be as in the statement; the given assertion
easily reduces to the following: there exists $s \in \integers$ such that
when we write $\gk$ as a linear combination of simple roots, the
coefficients of $\al_p, p \in  M_k(l-1,r-1)$ are all equal to $s$.

Let $s:=\dep(\gamma)$; if $\gk = \sum c_p \al_p$, remark \ref{staremk}
gives $c_u = c_v =s$. Now, all this holds for any $u, v$ which satisfy
$i(u) \geq \ls(\g)$, $\bi(v) \geq \rs(\g)$ and $i(v) = i(u) + 1$. This
concludes the proof. \qed

The figure now shows the nodes of $Z_k$ labeled by $b_i, c_j, s$:
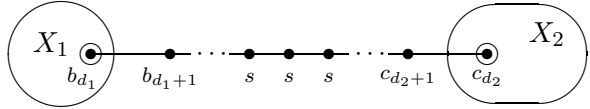
$$\begin{picture}(5,3)(3,-1)
\put(-1,0.25){\makebox(0,0){\bf \large $X_1$}}
\put(-.75,0){\circle{3}}
\put(0,0){\circle*{.25}}
\put(0,0){\circle{.5}}
\put(2,0){\circle*{.25}}
\put(4,0){\circle*{.25}}
\put(5,0){\circle*{.25}}
\put(6,0){\circle*{.25}}
\put(8,0){\circle*{.25}}
\put(10,0){\circle*{.25}}
\put(10,0){\circle{.5}}
\put(10.75,0){\oval(3.5,2.5)}
\put(11.5,0.5){\makebox(0,0){\bf \large $X_2$}}
\put(0,0){\line(1,0){2.5}}
\put(10,0){\line(-1,0){2.5}}
\put(3.5,0){\line(1,0){3}}
\put(2.6,0){\makebox(1,0){$\cdots$}}
\put(6.6,0){\makebox(1,0){$\cdots$}}
\put(-0.2,-.6){\makebox(0,0){\footnotesize $b_{d_1}$}}
\put(2,-.6){\makebox(0,0){\footnotesize $b_{d_1+1}$}}
\put(8,-.6){\makebox(0,0){\footnotesize $c_{d_2+1}$}}
\put(10,-.6){\makebox(0,0){\footnotesize $c_{d_2}$}}
\put(4,-.6){\makebox(0,0){\footnotesize $s$}}
\put(5,-.6){\makebox(0,0){\footnotesize $s$}}
\put(6,-.6){\makebox(0,0){\footnotesize  $s$}}
\end{picture}$$

\subsection{Main theorem}

Our main theorem is the following generalization of Theorem~(4.5) of
\kv.
\begin{theorem}\label{mainthm}
Let \xp be an extensible pair. 
\be
\item Given $\L, \mu, \nu \in \Htp$, there
exists an integer $K = K(\L, \mu, \nu)$ such that $\clmn =
\clmn[k\pr]$ for all $k, k\pr \geq K$.
\item Given $\L \in \Htp$ and $\beta \in \Ht$, there exists $K$ such
  that $\blb(k) = \blb(k\pr)$ for all $k, k\pr \geq K$.
\ee
\end{theorem}
{\bf Proof:}
Let $l$ be the maximum of the left supports of $\L, \mu, \nu, \beta$
and $r$ be the  maximum of their right supports.
Let $X_1\pr :=X_1(l-d_1)$ be the diagram formed by the nodes $p \in
L_k(l)$ and $X_2\pr:=X_2(r-d_2)$ be that formed by $p \in R_k(r)$.
 We can view
$\L, \mu, \nu, \beta$ as elements of $W(X_1\pr, X_2\pr)$ as in section
 \ref{twofour}. 

Define $\gamma:=\lambda + \mu - \nu$. Corollary \ref{corr} implies that
if $(\lambda  +\mu) \sim \nu$, then  $\nob=0$. Proposition \ref{aaa}
now applies and we deduce stabilization of $\clmn$ by applying
Theorem \ref{mainthm1} . 

The second assertion regarding the $\blb(k)$ follows analogously. We
set $\gamma = \lambda - \beta$ and use Proposition \ref{aaa} and
Theorem \ref{mainthm1}. \qed

\section{The stable representation ring}\label{stabrepring}
\subsection{Partial orders}
Having established that the multiplicities $\clmn$ stabilize for
extensible pairs $\xp$, we
shall now use the stable values $\clmn[\infty]$ as structure constants
to define an associative multiplication  $*$. First, we take a small
detour and define a partial order on the set $\Htp$.

Let \xp be an extensible pair of marked Dynkin diagrams as before.
The set $\pzk[+]$ of dominant weights of $\kma(Z_k)$ 
is partially ordered  with
$$\tau_1, \tau_2 \in \pzk[+], \;\; \tau_1 \geq \tau_2 \; \iff \; \tau_1 -
\tau_2 \in \qzk[+]$$

We can also make $\Htp$ into a poset. Let $\lambda_1 = (x,y), \, \lambda_2
= (z,w) \in \Htp$ be such that 
$\nob[\lambda_1] = \nob[\lambda_2]$. Let $l = \max(\ls(\lambda_1), \ls(\lambda_2))$
and $r = \max(\rs(\lambda_1), \rs(\lambda_2))$.
Proposition~\ref{aaa} implies that
there exist integers $b_i \; (1\leq i \leq l-1$) , $c_j \; (1\leq j \leq r-1)$
and $s$ such that for $k \geq l+r - d_1 - d_2$
\begin{equation}\label{biciagain}
\L_1\rk - \L_2\rk = \sum_{p \in L_k(l-1)} b_{i(p)} \al_p +  \sum_{p \in M_k(l-1,
  r-1)} s\, \al_p +  \sum_{p \in R_k(r-1)} c_{\bi(p)} \al_p
\eeq

\begin{definition}
Given $\lambda_1, \lambda_2 \in \Htp$, define
$\lambda_1 \po \lambda_2$ iff $\nob[\lambda_1]| = \nob[\lambda_2]$
and the
$b_i, c_j, s$ which occur in Equation~\eqref{biciagain} are all
non-negative.
\end{definition}

It is easy to check that $\po$ is a partial order on $\Htp$              
and that $\lambda_1 \po \lambda_2$ implies that 
$\lambda_1+\mu \po \lambda_2+\mu$ for all $\mu \in \Htp$.
We also have these equivalent
conditions which follow easily:
\begin{align} \label{eqcondgx}
\lambda_1 \po \lambda_2 &\Leftrightarrow \lambda_1\rk \geq \lambda_2\rk
\;\; \forall k \geq l+r -d_1 - d_2 \\
&\Leftrightarrow \lambda_1\rk \geq  \lambda_2\rk
 \; \text{ for infinitely many values of  } k \label{infpo}
\end{align}
where $\geq$ is the partial order on $\pzk[+]$.

\noindent
We also have:
\begin{lemma}\label{simplepropsofgx}
\begin{enumerate}
\item If $\L_1 \po \L_2$, then $\dep(\L_1 - \L_2) = s \geq 0$
\item If $\L, \mu, \nu \in \Htp$ are such that $\clmn[\infty] > 0$, then 
$\L + \mu \po \nu$.\\
\end{enumerate}
\end{lemma}
\noindent
{\bf Proof:} The proof of (1) is trivial. For (2), observe that
$\clmn[\infty] > 0$ implies that $\clmn >0$ for all large $k$. In
particular
$\L\rk + \mu\rk \geq \nu\rk$ for all large $k$. Equation~\eqref{infpo}
completes the proof. \qed
\subsection{One sided intervals}
The goal of this subsection is to analyze the partial orders 
on $\Htp$ and $\pzk[+]$.

Assume \xp is an extensible pair. 
For each $k \geq 1$, we know that $\kma(Z_k)$ is a
symmetrizable Kac-Moody algebra. Suppose $\xp \neq (A_{d_1}, A_{d_2})$
and $k$ is such that $\det Z_k \neq 0$, then
$Z_k$ cannot be a {\em finite} type \dyn; to see this we note
that among the series of finite type Dynkin diagrams, only
$A_n$ is extensible while the $B/C/D$ types are not (see Table 1 of
\kv). The $\det Z_k \neq 0$ condition also implies that $Z_k$ cannot be {\em
  affine}. So, it must be of {\em indefinite} type.

We first derive a result about the poset of dominant integral weights
in any Kac-Moody algebra of {\em indefinite} type.
Let $\kma$ be a symmetrizable Kac-Moody algebra of
 indefinite type. Let 
 $A$ be the corresponding $n \times n$ generalized Cartan
 matrix. Assume $A$ is indecomposable and that $\det A \neq 0$. We
 use the usual notation for roots,  weights etc. 
Let $(.|.)$ denote a nondegenerate,
 symmetric, bilinear, $\complex$-valued form on the Cartan subalgebra
$\csa$. 
This form exists because $\kma$ is symmetrizable.

Fix $\gamma \in P^+(\kma)$. Consider the set:
$$U(\gamma) := \{ \lambda \in P^+(\kma): \lambda \geq \gamma\}$$
Here again $\lambda \geq \gamma$ is defined by the condition
$\lambda - \gamma \in Q^+(\kma)$. 
The set $U(\gamma)$ is  the {\em one-sided interval} in the
poset $P^+(\kma)$, bounded below
 by $\gamma$.

\begin{proposition}\label{ugfinset}
$U(\gamma)$ is a finite set.
\end{proposition}

\noindent
{\bf Proof:} The definition of an {\em indefinite type} Kac-Moody
algebra \cite[Chapter 4]{Kac} implies that there exists $\xi =
\sum_{i=1}^n u_i \alpha_i \,\in \csa^*$ which satisfies the conditions\\
 (a) $u_j > 0 \,\forall j$ \\
 (b) $\xi(\Check{\alpha}_j) < 0 \, \forall j$.

\noindent
Condition (b) implies that $(\xi|\alpha_j) < 0\, \forall j$ since  
$(\xi|\alpha_j) = \frac{(\alpha_j|\alpha_j)}{2} \,\xi(\Check{\alpha}_j)$
and \\$(\alpha_j|\alpha_j) > 0$. 

\noindent
Condition (a) together with the fact that $\omega_j(\Check{\alpha}_i)
= \delta_{ij}$ implies that $(\xi|\omega_j) > 0\, \forall j$. To see
this observe:
\begin{align} 
(\xi|\omega_j) &= \sum_{i=1}^n u_i (\alpha_i|\omega_j) \notag\\
               &= u_j (\alpha_j|\alpha_j)/2 \notag
\end{align}
 
\noindent
For any $\eta \in \csa^*$, define $\T(\eta) := (\eta|\xi)$. The
conclusions of the above paragraph imply:\\
{\bf Fact 1}: $\T(\alpha) \leq 0$ for all $\alpha \in Q^+(\kma)$.\\
{\bf Fact 2}: $\T(\lambda) \geq 0$ for all $\lambda \in P^+(\kma)$.

\vspace{0.1in}
\noindent
If $\lambda \in U(\gamma)$, then $\lambda - \gamma \in Q^+(\kma)$. 
Hence $\T(\lambda) - \T(\gamma) = \T(\lambda - \gamma) \leq
0$. So, we have:\\
{\bf Fact 3}: If $\lambda \in U(\gamma)$, then 
$$0 \leq \T(\lambda) \leq \T(\gamma)$$

\noindent
We can now show that $U(\gamma)$ is finite. Let 
$\lambda = \sum_{i=1}^n a_i \omega_i$ with $a_i \in \N \,\forall
i$. Then $\T(\lambda) = \sum_{i=1}^n a_i (u_i (\alpha_i|\alpha_i)/2)
\leq \T(\gamma)$. Since each term in the sum is nonnegative, this
means that 
$$a_i \leq  \frac{2 \,\T(\gamma)}{u_i (\alpha_i|\alpha_i)} \; \forall i$$
Since $a_i$ is a nonnegative integer, there are only finitely
many choices for $a_i$. Thus $U(\gamma)$ is a finite set. \qed


\begin{remark}
The above Proposition is false if $\kma$ is of Finite or Affine
type. For instance if $\kma = \kma(A_1) = sl_2 \complex$, and $\gamma
= \alpha_1$, then $U(\gamma) = \{n\alpha_1: n \geq 1\}$. Similarly for
the rank 2 affine Lie algebra $\kma = A_1^{(1)}$ (affine $sl_2$), if
$\gamma = \alpha_1 + \alpha_2$, then $U(\gamma) = 
\{n(\alpha_1 + \alpha_2): n \geq 1\}$.
\end{remark}

\noindent
We now consider the poset $(\Htp,\po)$ and ask the analogous
question. Fix $\gamma \in \Htp$ and let 
\begin{equation}\label{udef}
U(\gamma) := \{ \lambda \in \Htp: \lambda \po \gamma\}
\end{equation}
We use the same notation $U(\cdot)$ that we used before, but the context
will resolve any confusion.
The following example shows that the poset $(\Htp, \po)$ is unlike
the poset $P^+(\kma)$; even if most (or all)
 members of the set
$\{\kma(Z_k)\}$ are of indefinite type, the set $U(\gamma)$ may still be infinite.

\begin{example}\label{eginft}
Let $X_1 = E_{10}$ with the node numbering and distinguished vertex
shown below
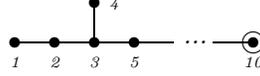
\setlength{\unitlength}{15pt}
$$\begin{picture}(4,3)(1,-1)
\put(0,0){\circle*{.25}}
\put(1,0){\circle*{.25}}
\put(2,0){\circle*{.25}} 
\put(2,1){\circle*{.25}}
\put(3,0){\circle*{.25}} 
\put(6,0){\circle*{.25}} 
\put(6,0){\circle{.5}} 
  \put(0,0){\line(1,0){4}}
  \put(2,0){\line(0,1){1}}
  \put(6,0){\line(-1,0){1}}
  \put(4.05,0){\makebox(1,0){...}}
  \put(0,-.5){\makebox(0,0){\tiny 1}}
  \put(1,-.5){\makebox(0,0){\tiny 2}}
  \put(2,-.5){\makebox(0,0){\tiny 3}}
  \put(2.5,1){\makebox(0,0){\tiny 4}}
  \put(3,-.5){\makebox(0,0){\tiny 5}}
  \put(6,-.5){\makebox(0,0){\tiny 10}}
  \end{picture}
$$
and take $X_2 = A_1$. We let $\epsilon_i$ denote the element 
$(0,0,\cdots,0, 1, 0,0 \cdots)$ of $\hop$, where the
$1$ occurs in the $i^{th}$ position. Take $\gamma :=(\epsilon_{10},0)$. 
For $k > 0$, consider the element $\lambda_k := (\epsilon_{10+k}, \epsilon_k)$
of $\Htp$. We claim that $\lambda_k \in U(\gamma)$ for all
$k > 0$. The proof follows by actual computation of the $b_i, c_j, s$:
\begin{align*}
b_i &= 0  \;\;\; (i \leq 10), & b_i &= i -10 \;\;\; (10 < i < 10+k)\\
s&=k & c_i &=  i \;\;\; (1 \leq i < k)
\end{align*}
Thus $U(\gamma)$ is an infinite set. \qed
\end{example}

However, a weaker finiteness assertion is true for the poset $\Htp$.
Given $s \in \integers$, define 
$$U(\gamma,s) := \{\lambda \in U(\gamma): \dep(\lambda - \gamma)
  =s\}$$
Lemma~\ref{simplepropsofgx} implies that $U(\gamma,s) = \emptyset$
  if $s<0$. Clearly
$$U(\gamma) = \bigsqcup_{s=0}^{\infty} U(\gamma,s)$$
We have:
\begin{proposition}\label{ugsfin}
$U(\gamma,s)$ is a finite set for each $s \in \N$.
\end{proposition}

\noindent
{\bf Proof:} Let $\gamma = (x,y)$, and suppose $\lambda =(z,w) \in
U(\gamma,s)$  with $x,y,z,w \in \hop$. Write
$x=(x_i)_{i \geq 1},\, y=(y_i)_{i \geq 1},\,
z=(z_i)_{i \geq 1},\, w=(w_i)_{i \geq 1}$. 
Let $l=\max(\ls(\gamma), \ls(\lambda))$ , $r = \max(\rs(\gamma),
\rs(\lambda))$.

First assume that  $\xp = (A_{d_1}, A_{d_2})$. Then $\D[1] = \D[2] =1 $ and
$\aip[]{1} = \aip[]{2} = i \; \forall i$. We have
$$s = \sum_j j z_j - \sum_j j x_j =  \sum_j j w_j - \sum_j j y_j $$
Clearly, there are  only finitely
many choices for $(z,w)$ such that this condition holds. Hence the number
of $\lambda = (z,w)$ in $U(\gamma,s)$ is finite and we are done.

So, assume now that $\xp \neq (A_{d_1}, A_{d_2})$.
Since $\lambda \po \gamma$, there are non-negative integers 
$b_i \; (1\leq i \leq l-1$) , $c_j \; (1\leq j \leq r-1)$
and $s$ such that for all large $k$
\begin{equation}\label{usualbici}
\L\rk - \g\rk =   \sum_{p \in L_k(l-1)} b_{i(p)} \al_p +  \sum_{p \in M_k(l-1,
  r-1)} s\, \al_p +  \sum_{p \in R_k(r-1)} c_{\bi(p)} \al_p
\end{equation}
\begin{claim}\label{themainclaim}
$\ls(\L) \leq \ls(\g) + s$ and $\;\;\rs(\L) \leq \rs(\g) + s$.
\end{claim}

\noindent
{\bf Proof:} Suppose $\ls(\L) > \ls(\g) + s$.
Then $l = \ls(\L)$. Take $p \in Z_k$ such that $i(p) =l$. Then
$(\L\rk - \g\rk)(\alch_p) > 0$.
From
Equation~\eqref{usualbici}, we get:
$$(\L\rk - \g\rk)(\alch_p) = s - b_{l-1}$$
Thus $$({\bf a})\;\; b_{l-1} < s$$
Again for all $q \in Z_k$ such that $\ls(\gamma) < i(q) < \ls(\L)$, we
have
$(\L\rk - \g\rk)(\alch_q) = \L\rk(\alch_q) \geq
0$. Using Equation~\eqref{usualbici} again, we get:
$$(\L\rk - \g\rk)(\alch_q) = 2b_{i(q)}-b_{i(q)+1}-b_{i(q)-1}$$
Hence 
$$({\bf b}) \;\;b_{i(q)-1} \leq 2b_{i(q)}-b_{i(q)+1}$$
Facts ({\bf a}) and
({\bf b}) imply that  $\{s,b_{l-1},b_{l-2},\cdots,
b_{\ls(\g)}\}$ is a strictly decreasing sequence of
nonnegative integers. But the number of terms in the sequence is \\
$l-\ls(\g)+1  > s+1$, which is a contradiction. The proof for the
right support is analogous. Our claim is thus true. \qed

So, the support $\supp(\L) \leq \supp(\g) + 2s$. Choose $k$ large  such that
(i) $\supp(\g) + 2s \leq d_1 + d_2 + k$ {\em and} 
(ii) $\kma(Z_k)$ is a Kac-Moody algebra of
indefinite type. This is possible, since we started with the
assumption that $\xp \neq (A_{d_1}, A_{d_2})$.

By our claim above, given any $\L \in U(\g, s)$, 
$\L\rk$ is a well defined dominant weight of
$\kma(Z_k)$. Further we also have $\L\rk \geq \g\rk$ in $\pzk[+]$,
i.e, $\L\rk \in U(\g\rk) \subset \pzk[+]$.
But Proposition~\ref{ugfinset} implies that the set  $U(\g\rk)$ is finite.
This implies that the number of possible choices for $\L$ is also
finite. \qed

\subsection{}
In this subsection, we will define a commutative (associative) $\complex$
algebra $\rxx$ whose multiplication operation has the $\clmn[\infty]$
as structure constants. 

\begin{definition}
If $\lambda = (x,y) \in \Ht$, we define
$$\bxo := \frac{\sum_{i} x_i \aip[]{1}}{\D[1]}$$ and
$$\bxt:= \frac{\sum_{i}  y_i\aip[]{2}}{\D[2]}$$
\end{definition}
Clearly $\nob[\lambda]= \D[1]\D[2] (\bxo - \bxt)$ and if
$\nob[\lambda] =0$, then $\dep(\lambda) = \bxo = \bxt$.

Let \mcl be the $\complex$ vector space with basis $\{v_\L : \L \in
\Htp\}$ and let \mclh denote its formal completion; so a typical
element of \mclh is an infinite sum $\xi:=\sum_{\L \in \Htp} c_\L v_\L$.
\begin{definition}
Let $\rxx$ be the set of all  $\xi:=\displaystyle\sum_{\L \in \Htp} c_\L v_\L \in
\mclh$ which satisfy the following two conditions:
\be
\item $\{\nob[\L] : c_\L \neq 0\}$  is a finite set
\item $\{\bxo : c_\L \neq 0\} \subset \frac{1}{\D[1]} \integers$  is bounded
above.
\ee
\end{definition}
\begin{remark}
 Since $$\bxt = \bxo - \frac{\nob[\L]}{\D[1]\D[2]}$$ the two
  conditions above imply that  $\{\bxt : c_\L \neq 0\} \subset
  \frac{1}{\D[2]} \integers$  is also bounded above.
\end{remark}

\subsection{The operation $*$ on $\rxx$}
{\bf STEP 1:} Observe that $v_\L, v_{\mu} \in \rxx$ for all $\L, \mu
\in \Htp$. Define
$$v_\L * v_\mu := \sum_{\nu \in \Htp} \clmn[\infty] \,v_\nu$$
The right hand side is in $\rxx$ since $\clmn[\infty] > 0$ implies $\L
+ \mu \po \nu$, i.e, 
\begin{align}
\nob[\nu] &= \nob[\L] + \nob[\mu] \text{  and } \label{clmone}\\
\dep(\L + \mu - \nu) &= \bxo + \bxo[\mu] - \bxo[\nu] \geq 0 \label{clmtwo}
\end{align}

\noindent
{\bf STEP 2:} Given  $\xi:=\sum_{\L \in \Htp} c_\L v_\L$ and
 $\eta:=\sum_{\mu \in \Htp} d_\mu v_\mu$ in $\rxx$, we define $\xi *
 \eta$ by bilinearity, i.e,
\begin{align}
\xi * \eta &:=\sum_{\lambda}\sum_{\mu} c_\L \,d_\mu (v_\L * v_\mu)\notag\\
&= \sum_\L \sum_\mu \sum_{\nu}  c_\L\, d_\mu\, \clmn[\infty]\,
v_\nu \label{rhson}
\end{align}
We need to  show that \eqref{rhson} is well defined. More precisely,
we need to show the following:
\begin{proposition}
Let  $\xi:=\sum_{\L \in \Htp} c_\L v_\L$ and
 $\eta:=\sum_{\mu \in \Htp} d_\mu v_\mu$ be elements of $\rxx$. Given
 $\nu \in \Htp$, the set $$F_\nu :=\{(\L, \mu) \in \Htp \times \Htp : 
c_\L \neq 0, d_\mu\neq 0, \clmn[\infty] \neq 0 \}$$ is finite.
\end{proposition}
\noindent
{\bf Proof:} As before, $\clmn[\infty] > 0$ gives us
Equations~\eqref{clmone}, \eqref{clmtwo}. Since $\xi, \eta \in \rxx$, the sets
$\{ \bxo : c_\L \neq 0\}$ and $\{\bxo[\mu]: d_\mu \neq 0\}$ are
bounded above. Together with Equation~\eqref{clmtwo} this implies that
these sets are also bounded below; so in fact both sets must be
finite. So $\dep(\L + \mu - \nu)  = \bxo + \bxo[\mu] - \bxo[\nu]$
takes only finitely many values $s_1, s_2, \cdots, s_p$ say.
This means that 
$$\L + \mu \in \bigcup_{i=1}^p U(\nu, s_i)$$
By Proposition~\ref{ugsfin} the right hand side is a finite
set. Thus $F_\nu$ is finite too. \qed

Thus * is a well defined operation on $\rxx$.

\begin{remark}
Let $X$ be an extensible diagram; definition~(5.5) of 
\kv introduced a $\complex$ algebra $\Lambda^{X}$ which is
 a subspace of  \mclh. From its definition, it  is clear that 
$\Lambda^{X}$ is just the subalgebra of $\rxxa$ generated by the
 $\{v_\L : \L \in \Htp\}$.
\end{remark}

\subsection{Associativity}
In order to show that $*$ is associative, we need to show that for all
$\L, \mu, \nu \in \Htp$ , $(v_\L * v_\mu) * v_\nu = v_\L * (v_\mu *
v_\nu)$. This was proved in \kv (see equation~(5.10) there)
for the $X(m)$ series. The proof in \kv uses the so called {\em
  Interval stabilization lemma}. We state the corresponding lemma
for our context:
\begin{lemma}
 Let  $\lambda_1, \lambda_2 \in
  \mathcal{H}_2^+$ with $\lambda_1 \po \lambda_2$. Let $I(\L_1, \L_2)
 := \{ \gamma
  \in \mathcal{H}_2^+ : \lambda_1 \po \gamma \po \lambda_2\}$ and 
$I\rk(\L_1, \L_2) := \{ \beta  \in P^+(Z_k) : \lambda_1\rk \geq \beta \geq
  \lambda_2\rk\}$ for $k$ large. Then
\begin{enumerate}
\item $I(\L_1, \L_2)$ is a finite set
\item There exists $N$ such that for all $k \geq N$, $I\rk(\L_1, \L_2) = \{
  \gamma\rk: \gamma \in I(\L_1, \L_2) \}$
\end{enumerate}
\end{lemma}
The proof is exactly the same as the proof of Lemma~(5.1) of
\kv; instead of equation~(5.1) of \kv, we use our 
equation~\eqref{biciagain}. \qed

Now, given $\lambda, \mu, \nu \in \Htp$, we fix $\gamma \in \Htp$ and
let $\clmng$ denote the multiplicity of the representation $L(\gk)$ in
the triple tensor product $L(\lambda\rk) \otimes L(\mu\rk) \otimes
L(\nu\rk)$. If this multiplicity becomes a constant for large $k$, we
denote the constant value by $\clmng[\infty]$. It is easily seen that
the associativity of $*$ is implied by the following lemma:
\begin{lemma}
The multiplicity $\clmng$ does become constant for large $k$ and we
have
$$\clmng[\infty] = \sum_{\delta \in \Htp} c_{\lambda\mu}^{\,\delta}(\infty)
\, c_{\delta\nu}^{\,\gamma}(\infty) =  
\sum_{\delta \in \Htp} c_{\lambda\delta}^{\,\gamma}(\infty)\,
 c_{\mu\nu}^{\,\delta}(\infty)$$
\end{lemma}
The above relation is just equation~(5.8) of \kv. The proof given
there carries over with no change. \qed

We call $\rxx$ with the operation $*$ the {\em stable representation
  ring} for the pair $\xp$.

\section{BCD diagrams}\label{BCD}
Having dealt with extensible pairs in sections \ref{extpair} -
\ref{stabrepring}, we now turn to the classical $B_n, C_n, D_n$
sequences. These fall outside the class of extensible sequences. Our
goal is to give simple characterizations of weights $\L, \mu, \nu,
\beta$ which are of the form of Theorem \ref{mainthm1} and hence
exhibit stabilization behavior.

Let $X$ be one of  the diagrams $B_3, C_3, D_3$ below :
\be
\item $B_3$: 
\begin{picture}(2,1)(-1,0)
\put(0,0.5){\circle*{.25}}
\put(1.5,0.5){\circle*{.25}}
\put(2.5,0.5){\circle*{.25}}
\put(2.5,0.5){\circle{.5}}
\put(0,0.4){\line(1,0){1.5}}
\put(0,0.6){\line(1,0){1.5}}
\put(1.5,0.5){\line(1,0){1}}
\put(0.75,0.5){\makebox(0,0){$<$}}
\put(0,0){\makebox(0,0){\tiny 1}}
\put(1.5,0){\makebox(0,0){\tiny 2}}
\put(2.5,0){\makebox(0,0){\tiny 3}}
\end{picture}

\item $C_3$: 
\begin{picture}(2,1)(-1,0)
\put(0,0.5){\circle*{.25}}
\put(1.5,0.5){\circle*{.25}}
\put(2.5,0.5){\circle*{.25}}
\put(2.5,0.5){\circle{.5}}
\put(0,0.4){\line(1,0){1.5}}
\put(0,0.6){\line(1,0){1.5}}
\put(1.5,0.5){\line(1,0){1}}
\put(0.75,0.5){\makebox(0,0){$>$}}
\put(0,0){\makebox(0,0){\tiny 1}}
\put(1.5,0){\makebox(0,0){\tiny 2}}
\put(2.5,0){\makebox(0,0){\tiny 3}}
\end{picture}

\item $D_3$:
\setlength{\unitlength}{15pt}
\begin{picture}(3,1)(-1,0)
\put(0,0.5){\circle*{.25}}
\put(1,0){\circle*{.25}}
\put(1,0){\circle{.5}}
\put(0,-0.5){\circle*{.25}}
\put(0,.5){\line(2,-1){1}}
\put(0,-.5){\line(2,1){1}}
\put(-.5,.5){\makebox(0,0){\tiny 1}}
\put(-.5,-.5){\makebox(0,0){\tiny 2}}
\put(1.5,0){\makebox(0,0){\tiny 3}}
 \end{picture}
\end{enumerate}

The nodes marked by the extra circles are taken to be distinguished.
Then $X(m) = Y_{m+3}$ where $Y \in \{B,C,D\}$ (see figure \ref{figsbcd}).

\begin{figure}
\begin{center}
\begin{picture}(6,1)(-1,0)
\put(0,0.5){\circle*{.25}}
\put(1.5,0.5){\circle*{.25}}
\put(2.5,0.5){\circle*{.25}}
\put(3.5,0.5){\circle*{.25}}
\put(5.5,0.5){\circle*{.25}}
\put(6.5,0.5){\circle*{.25}}
\put(0,0.4){\line(1,0){1.5}}
\put(0,0.6){\line(1,0){1.5}}
\put(1.5,0.5){\line(1,0){2.5}}
\put(6.5,0.5){\line(-1,0){1}}
\put(0.75,0.5){\makebox(0,0){$<$}}
\put(4.75,0.5){\makebox(0,0){$\cdots$}}
\put(0,0){\makebox(0,0){\tiny 1}}
\put(1.5,0){\makebox(0,0){\tiny 2}}
\put(2.5,0){\makebox(0,0){\tiny 3}}
\put(3.5,0){\makebox(0,0){\tiny 4}}
\put(6.5,0){\makebox(0,0){\tiny $n$}}
\end{picture}\\
\begin{picture}(6,1.5)(-1,0)
\put(0,0.5){\circle*{.25}}
\put(1.5,0.5){\circle*{.25}}
\put(2.5,0.5){\circle*{.25}}
\put(3.5,0.5){\circle*{.25}}
\put(5.5,0.5){\circle*{.25}}
\put(6.5,0.5){\circle*{.25}}
\put(0,0.4){\line(1,0){1.5}}
\put(0,0.6){\line(1,0){1.5}}
\put(1.5,0.5){\line(1,0){2.5}}
\put(6.5,0.5){\line(-1,0){1}}
\put(0.75,0.5){\makebox(0,0){$>$}}
\put(4.75,0.5){\makebox(0,0){$\cdots$}}
\put(0,0){\makebox(0,0){\tiny 1}}
\put(1.5,0){\makebox(0,0){\tiny 2}}
\put(2.5,0){\makebox(0,0){\tiny 3}}
\put(3.5,0){\makebox(0,0){\tiny 4}}
\put(6.5,0){\makebox(0,0){\tiny $n$}}
\end{picture}\\
\begin{picture}(7,2)(-2,-1)
\put(0,0.5){\circle*{.25}}
\put(1,0){\circle*{.25}}
\put(2,0){\circle*{.25}}
\put(3,0){\circle*{.25}}
\put(5,0){\circle*{.25}}
\put(6,0){\circle*{.25}}
\put(0,-0.5){\circle*{.25}}
\put(0,.5){\line(2,-1){1}}
\put(0,-.5){\line(2,1){1}}
\put(1,0){\line(1,0){2.45}}
\put(6,0){\line(-1,0){1.45}}

\put(4.1,0){\makebox(0,0){$\cdots$}}

\put(-.5,.5){\makebox(0,0){\tiny 1}}
\put(-.5,-.5){\makebox(0,0){\tiny 2}}
\put(1.25,-0.5){\makebox(0,0){\tiny 3}}
\put(2,-0.5){\makebox(0,0){\tiny 4}}
\put(6,-0.5){\makebox(0,0){\tiny $n$}}
 \end{picture}
\caption{$B_n, C_n, D_n$}\label{figsbcd}
\end{center}
\end{figure}

It is an easy fact that $\det(B_n) = 2 = \det(C_n)$ and $\det(D_n)=4$
for all $n \geq 3$. So $\Delta_X=0$ and $X$ is not an extensible
diagram. In fact the proofs of the results of \kv often 
used the fact that for extensible
$X$,   $\Delta_X \neq 0$ and hence $|\det X(m)| \to \infty$ as $m \to
\infty$.

Recall the definitions of the sets $\ho, \hop, \Ht, \Htp$ etc from section
\ref{notats}. For $Y \in \{B,C,D\}$, we let $\omi$ (resp  $\ali$)
be the fundamental  weight (resp. simple root) corresponding to the $i^{th}$
node of $Y_n$ in figure \ref{figsbcd}.
We let $\ombi := \omi[n-i+1]$ and $\albi := \ali[n-i+1]$.
Given $\gamma = (x,y) \in \Ht$ and $n \geq \ell(x) + \ell(y)$, let 
$$ \gamma\rk[n] := \sum_{i} x_i \,\omi + \sum_i y_i \,\ombi $$
For the series $B_n, C_n, D_n$,  our aim is to characterize $\L, \mu,
\nu \in \Htp$, $\beta \in \Ht$  which have the specific form of
Theorem~\ref{mainthm1}.

Let $\gamma = (x,y) \in \Ht$; set $\noba := \displaystyle\sum_i
i\,y_i$. We let $l:= \max\{ i: x_i \neq 0\}$ and 
 $r:= \max\{ i: y_i \neq 0\}$.

\subsection{$B_n$}

First, we state  the following lemma concerning the fundamental weights of $B_n$
\begin{lemma}\label{lem1}
\be
\item Let $1 \leq i < n$. If $\ombi = \displaystyle\sum_{j=1}^n 
 c_{ij}(n) \, \albi[j]$, then $c_{ij}(n) = i$ for $i \leq j \leq n$.
\item $2 \omi[1] = \displaystyle\sum_{j=1}^{n} j\, \albi[j]$.
\item Let $1 < i \leq n$. If $\omi = \displaystyle\sum_{j=1}^n
 d_{ij}(n) \, \albi[j]$, then
$d_{ij}(n) = j$ for $1 \leq j \leq (n-i+1)$.
\ee
\end{lemma}

\noindent
{\bf Proof:} Follows by direct computation. \qed

For $\gamma$ as above,  define
$$\hht[B](\gamma):=\frac{x_1}{2} + \sum_{i>1} x_i$$

Suppose we write
$$ \gk[n] = \sum_{i=1}^n p_i(n) \, \albi \text{ with } p_i(n) \in \integers,$$
lemma \ref{lem1} implies 
\begin{equation}\label{eqabove}
p_i(n) = \sum_j j\,y_j + i\,\hht[B](\gamma) \;\; \text{  for } 
r \leq i \leq (n-l+1)
\end{equation}
We note that $\gk[n] \in Q^+(B_n)$ for all large $n$ implies that $\hht[B](\gamma)
\in \integers^{\geq 0}$.
Observe that $p_i(n)$ grows linearly as a function of $i$ unless
$\hht[B](\gamma)=0$, in which case 
$$p_i(n) = p_j(n) = \noba \;\;\; \text{ for } r \leq i,j \leq
(n-l+1)$$
This is exactly the requirement of Theorem \ref{mainthm1}. To rephrase
this in the notations of section \ref{twofour}, identify $B_l$ and $A_r$
with the subdiagrams of $B_n$ formed by the leftmost $l$ and the
rightmost $r$ nodes (see figure \ref{figsbcd}). Then equation \eqref{eqabove} shows
that 
$$\exists\, s \in \integers \text{ such that } \gamma \in R_s(B_l, A_r)
\Leftrightarrow \hht[B](\gamma)=0$$
In this case $s=\noba$. Theorem \ref{mainthm} now implies:
\begin{proposition}\label{bstab}
Consider the sequence $B_n$; let $\L, \mu, \nu \in \Htp$, $\beta \in \Ht$. 
\be
\item If $\hht[B](\L) + \hht[B](\mu) =  \hht[B](\nu)$, $\clmn[n]$
  stabilizes.
\item If $\hht[B](\L) = \hht[B](\beta)$, $\blb(n)$
  stabilizes.  
\ee
\end{proposition}
\begin{remark}
From the remarks of the above paragraph we conclude that if $\clmn[n]
> 0$ (resp. $\blb(n)>0$) for all large $n$, then $\hht[B](\L) +
\hht[B](\mu) \geq  \hht[B](\nu)$ (resp.  $\hht[B](\L) \geq
\hht[B](\beta)$). Thus our theorem above deals with the $\nu$ and
$\beta$ which have the maximum allowed height, and shows that
stabilization holds in this case.
\end{remark}

\begin{corollary}\label{onebcd}
Let $\L$ be a weight supported on the $A_r$ portion of the diagram of
$B_n$, i.e, let $\L =({\bf{0}}, y)$ for some $y \in \hop$. Then
$\hht[B](\L)=0$. If $\mu, \nu, \beta$ are also of this form, then the
height compatibility conditions of the above proposition are trivially
satisfied. So, for all such weights of the $B_n$'s, tensor product and
branching multiplicities stabilize.
\end{corollary}

\begin{example}
For this example, we
will use the following indexing scheme: let
 $[m_1 \,m_2\, \cdots\, m_n]$ denote the irreducible 
representation of $B_n$ with highest weight $\sum_{i=1}^n
m_i \omi$. 

Consider the representation $[100\cdots01]$ of $B_n$. Its highest
weight is obtained as $\Lk[n]$ where $\L:=(x,y)$ with $x_1 = y_1 :=1$
and $x_i = y_i = 0 \; \forall i > 1$. Observe that $\hht[B](\L) =
1/2$. 
We give below the
decomposition of the tensor square of this representation in $B_3, B_4$ and $B_5$. 
Note that our theorem above guarantees stabilization of multiplicities
for $\nu$ for which $\hht[B](\nu) = 1/2 + 1/2 = 1$. The data below was
generated using the program LiE:
\be
\item In $B_3$:
\begin{align}
[101]\otimes[101] &= 1.[202] +
1.[210] +
2.[201] +
2.[200]  \notag\\
& +1.[012] +
2.[011] +
2.[010] \\
&+ 1.[001] 
+1.[020] +
1.[003] +
1.[002] +
1.[000]  \notag
\end{align}
\item In $B_4$:
\begin{align}
[1001]\otimes[1001] &= 
1.[2002] +
1.[2010]+
2.[2001] +
2.[2000]\notag\\
&+1.[0102] + 
1.[0110] +
2.[0101] +
2.[0100] \\
&+1.[0012] +
2.[0011] +
2.[0010] \notag\\
&+1.[0020] +  
1.[0003] +
1.[0002]+
1.[0001] +
1.[0000] \notag
\end{align}
\item In $B_5$:
\begin{align}
[10001]\otimes[10001] &= 
1.[20002] +
1.[20010] +
2.[20001] +
2.[20000] \notag\\
&+1.[01002] +
1.[01010] +
2.[01001] +
2.[01000] \\
&+1.[00102] +
2.[00101] +
2.[00100] +
1.[00110] \notag\\
&\negthickspace\negthickspace
\negthickspace\negthickspace\negthickspace\negthickspace\negthickspace\negthickspace
\negthickspace\negthickspace\negthickspace\negthickspace\negthickspace\negthickspace
\negthickspace\negthickspace\negthickspace\negthickspace\negthickspace\negthickspace
\negthickspace\negthickspace\negthickspace\negthickspace\negthickspace\negthickspace
+2.[00011] +
1.[00002] +
2.[00010] +
1.[00012] +
1.[00020] +
1.[00003] +
1.[00001] +
1.[00000]\notag
\end{align}
\ee

In each case, the right hand side is arranged so that the weights
$\nu$ in the first three rows satisfy $\hht[B](\nu) =1$. Suppose we look at
such $\nu \in \Htp$ which make sense in $B_3, B_4$ and $B_5$; from the data
we do see that the multiplicities of these  are the same for
$B_n$, $n=3,4,5$. The last rows consists of $\nu$ for which $\hht[B](\nu) = 0$.

\end{example}

\subsection{$C_n$, $D_n$}
Define $\hht[C](\gamma) = \sum_i x_i$. The analogous statements for
$C_n$ are contained in the following:
\begin{proposition}\label{cstab}
\be
\item $\exists\, s \in \integers$ such that $\gamma \in R_s(C_l, A_r)
\Leftrightarrow \hht[C](\gamma)=0$, in which
case $s=\noba$. 
\item Given $\L, \mu, \nu \in \Htp$, 
$\hht[C](\L) + \hht[C](\mu) =  \hht[C](\nu)\Rightarrow \clmn[n]$
  stabilizes.
\item Given $\beta \in \Ht$ with $\hht[C](\L) = \hht[C](\beta)$, then 
$\blb(n)$  stabilizes.  
\ee
\end{proposition}
\noindent
{\bf Proof:} This again follows from the expression for 
the fundamental weights $\omi$ and $\ombi$ of 
$C_n$ in the basis of simple roots. \qed

For type $D$, the height function turns out to be :
$$\hht[D](\gamma) : = \frac{x_1}{2} +  \frac{x_2}{2} + \sum_{i>2} x_i$$
The three assertions of proposition \ref{cstab} are true for $D_n$ if we
replace $C_l$ by $D_l$ and $\hht[C]$ by $\hht[D]$.
We leave the remaining details to the reader.

\begin{remark}
\be 
\item For {\em one-sided} weights $\L, \mu, \nu, \beta$ as in
  corollary \ref{onebcd}, we evidently obtain stabilization in types $C$ and $D$
  as well.
\item We saw above that the $p_i(n)$ in general grows linearly with $i$. Our
proof of stabilization however only works when the $p_i(n)$'s are
constant. Hence we only deduce stabilization under the height
compatibility condition of propositions \ref{bstab} and \ref{cstab}.
\item However, 
from computer generated data of multiplicities in tensor products,
it appears that stabilization still holds, even when the $p_i(n)$'s
are non constant in the middle. Our method of proof however fails in
this situation. 
\ee
\end{remark}

The last remark leads us to the following conjecture:

\vspace{.2cm}
\noindent
{\bf Conjecture:} For types $B,C,D$, $\clmn[n]$ stabilizes for
{\em all} triples $\L, \mu, \nu \in \Htp$.

\vspace{.2cm}
It seems likely that a similar statement holds for the branching
multiplicities as well.

\subsection{Final remarks}\label{finalrem}
We briefly mention two other papers that are related to the work of
the present article. 

In \cite{bklms}, Benkart, Kang, Lee, Misra and Shin considered
the affine diagrams $A_n^{(1)}$ and looked at  integral
weights of a fixed level that are parametrized by $\Ht  \times \integers$ (the extra
degree of freedom comes from the fact that the Cartan subalgebra of
$A_n^{(1)}$ is $n+2$ dimensional, and so a weight is not uniquely
determined by specifying its values on the simple coroots
alone). Given a level zero element $\gamma \in \Ht \times \integers$, they
worked out the condition that $\gamma$ must satisfy so that
$\gamma\rk[n] \in Q(A_n^{(1)})$ for all large $n$. They show that such
$\gamma$'s must have the form as in Theorem~\ref{mainthm1} (see
proposition 1.22 of \cite{bklms}).

Our proof of stabilization then applies (even though $A_n^{(1)}$ is
not strictly of the form $Z_k(X_1, X_2)$ for any $X_1, X_2$). So for
$A_n^{(1)}$ we obtain stabilization of 
$\clmn[n]$ and $\blb(n)$ for all $\L, \mu, \nu, \beta$,
subject to a compatibility condition on their levels.

In the present work, we have not mentioned an important
representation theoretic statistic - weight multiplicities.
The main objective of \cite{bklms} (and its follow up paper
\cite{bkls}) was to analyze the behavior of these numbers as $k \to
\infty$. In our notation,
given $\L \in \Htp$ and $\beta \in \Ht$, let $m_{Z_k}(\L, \beta)$ be
the dimension of the weight space with weight $\beta\rk$ in the
representation $L(\Lk)$ of $\kma(Z_k)$. These papers show that for the
case when $Z_k$ is either $A_n$ or $A_n^{(1)}$, $m_{Z_k}(\L, \beta)$
is a polynomial in $k$. They also consider other classical and affine
sequences and establish this polynomiality under the assumption that the
weights $\L$ and $\beta$ are ``one-sided'' in the sense of corollary \ref{onebcd}.

As mentioned in the introduction, Webster \cite{ben} has recently proved more general versions of 
 the stabilization results of the present article using quiver varieties and their connections to
representation theory. He has also proved the polynomiality of weight
 multiplicities for more general $Z_k$'s, thereby extending the result
 of \cite{bkls}. In fact \cite{ben} uses the result of \cite{bkls} for the $A_k$
 to prove it for all $Z_k$'s. However the method originally used in
 \cite{bkls} is very different in flavor than that in \cite{ben}.

\bibliographystyle{amsplain}
\bibliography{master-bibliography}
\nocite{*}
\end{document}